\documentclass[smallextended,referee,envcountsect]{svjour3}
\smartqed
\usepackage{graphicx}
\journalname{ }
\usepackage[latin1]{inputenc}
\usepackage{amsfonts}
\usepackage{amsmath}
\usepackage{cases}
\usepackage[ruled,vlined]{algorithm2e}
\usepackage{booktabs}
\usepackage{epstopdf}
\usepackage{textcomp}
\usepackage{hyperref}
\usepackage{subfigure}
\usepackage{multirow}
\usepackage{rotating}
\usepackage{lscape}
\newtheorem{assumption}{Assumption}[section]
\usepackage{color}

\begin{document}

\title{Anderson Acceleration of Derivative-Free Projection Methods for Constrained Monotone Nonlinear Equations}

\titlerunning{Anderson Acceleration of Derivative-Free Projection Methods}

\author{Jiachen Jin \and
	        Hongxia Wang \and
	        Kangkang Deng
}

\institute{Jiachen Jin \at
	jinjiachen@nudt.edu.cn
	              \and
	              Hongxia Wang, Corresponding author \at
	              wanghongxia@nudt.edu.cn
	           \and
	           Kangkang Deng \at
	           freedeng1208@gmail.com\\
	\\
	College of Science, National University of Defense Technology, Changsha, 410073, Hunan, China.
}

\date{Received: date / Accepted: date}
%The correct dates will be entered by the editor.

\maketitle

\begin{abstract}
The derivative-free projection method (DFPM) is an efficient algorithm for solving monotone nonlinear equations. As problems grow larger, there is a strong demand for speeding up the convergence of DFPM. This paper considers the application of Anderson acceleration (AA) to DFPM for constrained monotone nonlinear equations. By employing a nonstationary relaxation parameter and interleaving with slight modifications in each iteration, a globally convergent variant of AA for DFPM named as AA-DFPM is proposed. Further, the linear convergence rate is proved under some mild assumptions. Experiments on both mathematical examples and a real-world application show encouraging results of AA-DFPM and confirm the suitability of AA for accelerating DFPM in solving optimization problems.
\end{abstract}

%All acknowledgements should be placed in the back of the paper after Conclusions..

\section{Introduction}\label{Intro}
In this paper, we focus on solving the following monotone nonlinear equations with convex constraint:
\begin{equation}\label{P}
	F(x)=0,\ x\in \mathcal{C},
\end{equation}
where $F : \mathbb{R}^n \to \mathbb{R}^n$ is a continuous and monotone mapping, and $\mathcal{C}\subset \mathbb{R}^{n}$ is a closed convex set. The monotonicity of the mapping $F$ means that
\[
(F(x)-F(y))^{\top}(x-y)\geq 0,~\forall~x,y\in \mathbb{R}^n.
\]
The systems of nonlinear equations have numerous applications, such as chemical equilibrium systems \cite{MM87}, split feasibility problems \cite{QX05} and neural networks \cite{CZ14}. Further, some concrete application models in real world are monotone. For instance, compressed sensing is firstly formulated for a convex quadratic programming, and then for an equivalent monotone nonlinear equations \cite{XWH11}. Regularized decentralized logistic regression also can be expressed as monotone nonlinear equations \cite{JYTH22}. As observation techniques advance, observed date size expands, and the requirement for resolution of results increases, the scale of nonlinear equations enlarges accordingly.

Various iterative methods for solving \eqref{P} include Newton method \cite{SWQ02}, trust-region algorithm \cite{QTL04}, Levenberg-Marquardt method \cite{YJM23}, etc. Although these methods perform well theoretically and numerically, they have difficulties in dealing with large-scale equations due to the computation of Jacobian matrix or its approximation. In contrast, on the basis of the hyperplane projection technique for monotone equations \cite{SS98} and the first-order optimization methods for unconstrained optimization, many derivative-free projection methods (DFPM) have sprung up \cite{LLXWT22,MJJYH23,WSLZ23,YJJLW21} for convex-constrained monotone nonlinear equations, whose computational cost in each iteration is only to calculate function values.

The search direction and line search procedure are crucial for DFPM, and different constructions of theirs correspond to different variants of DFPM. Benefiting from the simple structure and low storage capacity of conjugate gradient methods (CGM), the conjugate gradient projection methods (CGPM), which are based on the design of the search direction in CGM, provide a class of competitive algorithms, for instance, CGPM \cite{LF19}, spectral CGPM \cite{IKRPA22}, three-term CGPM \cite{WA22}.
Meanwhile, different line search procedures may obtain different convergence properties. Some line search procedures have been proposed for DFPM in solving constrained monotone nonlinear equations (see \cite{AK15,LL11,OL18,ZZ06} for instance).
Although there have been many studies on the DFPM for solving problem \eqref{P}, almost all of these existing studies focus on specific algorithms. Only a few papers have discussed unified studies on this class of methods partially (see \cite{GM23,OL23}), which motivates us to center on a general framework of DFPM and its convergence analysis.

In order to construct more efficient numerical algorithms, a promising strategy that has recently emerged in a number of fields is to embed acceleration techniques in the underlying algorithms. Anderson acceleration (AA) was originally designed for integral equations \cite{Anderson65} and is now a very popular acceleration method for fixed-point schemes. AA can be viewed as an extension of the momentum methods, such as inertial acceleration \cite{Alvarez00} and Nesterov acceleration \cite{Nesterov83}. The idea differs from theirs in maintaining information of previous steps rather than just two last iterates, and update iteration as a linear combination of the information with dynamic weights.
Some studies have explored the connection between AA and other classical methods, which also facilitates the understanding of AA. For linear problems, Walker and Ni \cite{WN11} showed that AA is related to the well-known generalized minimal residual algorithm (GMRES \cite{SS86}). Potra and Engler \cite{PE13} demonstrated the equivalence between GMRES and AA with any mixing parameters under full-memory (i.e., $m = \infty$ in Algorithm \ref{algo2}). For nonlinear cases, AA is also closely related to the nonlinear GMRES \cite{WHS21}. Fang and Saad \cite{FS09} identified the relationship between AA and the multi-secant quasi-Newton methods.

Although AA often exhibits superior numerical performance in speeding up fixed-point computations with countless applications, such as reinforcement learning \cite{WCHZZQ23}, numerical methods for PDE \cite{SRX19} and seismic inversion \cite{Yang21}, it is known to only converge locally in theory \cite{MJ20,TK15}. The convergence analysis of most, if not all, existing methods require the involved function is continuously differentiable \cite{SRX19,TK15,WN11}. New results in Bian and Chen \cite{BC22} proved that AA for $m=1$ is Q-linear convergent with a smaller Q-factor than existing Q-factors for a class of nonsmooth fixed-point problem.
Moreover, they proposed a modified AA for the nonsmooth fixed-point problem based on the smoothing approximation, and proved that it owns the same R-linear convergence rate as the classical AA for continuously differentiable case. More recent results in Garner et al. \cite{GLZ23} proved that AA improves the R-linear convergence factor over fixed-point iteration when the operator is linear and symmetric or is nonlinear but has a symmetric Jacobian at the solution. Rebholz and Xiao \cite{RX23} investigated the effect of AA on superlinear and sublinear convergence of various fixed-point iteration, with the operator satisfying certain properties.

The efficient procedure of AA in solving wide applications further motivates us to investigate this technique to DFPM for solving problem \eqref{P}. Our main goal in this paper is hence to provide a globally convergent AA of general DFPM without any further assumptions other than monotonicity. Clearly, the work is an extension of recent inertial DFPM in \cite{IKRPA22,MJJYH23,WSLZ23} due to utilizing more information than just last two iterates. The main contributions of this paper are outlined below:

$\bullet$ An accelerated version of DFPM combined with AA (AA-DFPM) is proposed to solve convex-constrained monotone nonlinear equations. Fully exploiting the optimization problem structure, several modifications are added to the acceleration algorithm. To the best of our knowledge, this is the first application of AA in DFPM.

$\bullet$ A self-contained proof for the global convergence of AA-DFPM is given with no additional assumptions apart from monotonicity on the nonlinear mapping. We further discuss the convergence rate under some standard assumptions.

$\bullet$ The numerical experiments on large-scale constrained nonlinear equations and decentralized logistic regression demonstrate that AA-DFPM outperforms the corresponding DFPM in terms of efficiency and robustness.

The paper is organized as follows. In Section \ref{Alg}, we start by outlining the unified algorithmic framework of DFPM and its convergence results. Based on DFPM with the convergence gained, in Section \ref{AAA}, we further introduce AA and extend the acceleration technique to DFPM. The convergence analysis of AA-DFPM is established under some mild assumptions in Section \ref{Cov}. We report the numerical results of AA-DFPM on large-scale constrained nonlinear equations and a machine learning problem in Section \ref{NE}. The conclusion is given in Section \ref{Con}.

\noindent{\bf Notation.} Throughout the paper, we denote by $\|\cdot\|$ be the Euclidean norm on $\mathbb{R}^n$ and $F_k:=F(x_k)$. For a closed convex set $\mathcal{C}$, ${\rm dist}(x_k,\mathcal{C})$ denotes the distance from an iterate $x_k$ to $\mathcal{C}$ and the projection operator
$P_{\mathcal{C}}[x]={\rm argmin}\{\|z-x\|~|~z\in\mathcal{C}\}$. Furthermore, it has the nonexpansive property:
\[
\|P_\mathcal{C}[x]-P_\mathcal{C}[y]\|\leq \|x-y\|,~\forall~x,y\in \mathbb{R}^n.
\]

\section{Derivative-Free Projection Method}\label{Alg}
In this section, we review a comprehensive framework of DFPM and recollect its theoretical results. Throughout the paper, we assume that the solution set $\mathcal{S}$ of problem (\ref{P}) is nonempty.
\subsection{General Framework of DFPM}
The core of DFPM is the hyperplane projection technique \cite{SS98}. It projects the current iterate onto a hyperplane constructed based on the monotonicity of the mapping, which separates the current iterate from the solution effectively.
In general, for a given current iterate $x_k$, a search direction $d_k$ is computed first, then a stepsize ${\alpha _k}$ is calculated by a line search to satisfy
\[F(z_k)^\top(x_k-z_k) > 0,\]
where $z_k = x_k+\alpha_k d_k$. By the monotonicity of $F$, we have
\begin{equation}\label{mon}
	F(z_k)^\top(x^*-z_k)=(F(z_k)-F(x^*))^{\top}(x^*-z_k)\leq 0,~\forall x^*\in \mathcal{S}.
\end{equation}

Thus the hyperplane
\[H_k:=\{x\in\mathbb{R}^n|F(z_k)^\top(x-z_k) = 0\}\]
strictly separates the current iterate $x_k$ from any solution $x^*$. Projecting $x_k$ first onto the separating hyperplane $H_k$ then onto the feasible set $\mathcal{C}$, $x_{k + 1} = P_\mathcal{C}\left[P_{H_k}\left[ x_k\right]\right]$. Separation arguments show that ${{\rm dist}(x_{k},\mathcal S)}$ decreases monotonically with the increase of $k$, which essentially ensures the global convergence of DFPM.

From the above process, the determination of $d_k$ and $\alpha_k$ plays a crucial role in DFPM. Different choices of direction or stepsize lead to different variants of DFPM. As mentioned earlier, competitive DFPM includes CGPM \cite{LF19}, spectral CGPM \cite{IKRPA22} and three-term CGPM \cite{WA22}. We concentrate on a unified framework for DFPM in Algorithm \ref{algo1}.

\begin{algorithm}[H]
	\caption{General framework of DFPM}\label{algo1}
	\KwIn{initial point $x_0\in\mathcal{C}$, parameters $\gamma,~\sigma,~\epsilon>0$, $0<s_1\leq s_2$, $\rho\in(0,1)$, $0\leq t_1\ll t_2$, $\zeta\in(0,2)$. Set $k: = 0$.}
	
	{\bf Step 1.} Compute $F_k$. If $\|F_k\|<\epsilon $, stop. Otherwise, go to Step 2.
	
	{\bf Step 2.} Compute the search direction $d_k$ such that
	\begin{equation}\label{SDC}
		F_k^\top d_k\leq -s_1 \|F_k\|^2,
	\end{equation}
	\begin{equation}\label{TRP}
		\|d_k\| \leq s_2 \|F_k\|.
	\end{equation}
	
	{\bf Step 3.} Set $z_k = x_k+\alpha_k d_k$, where ${\alpha _k} = \gamma \rho ^{i_k}$ with $i_{k}$ being the smallest nonnegative integer $i$ such that
	\begin{equation}\label{ak}
		- F(x_k + \gamma \rho^i d_k)^\top{d_k} \geq \sigma \gamma \rho ^i P_{[t_1,t_2]}[\|F(x_k + \gamma \rho ^i d_k)\|]\| d_k \|^2.
	\end{equation}
	
	{\bf Step 4.} Yield the next iteration by
	\begin{equation}\label{xk1}
		x_{k + 1} = P_\mathcal{C}\left[ x_k - \zeta u_k F(z_k)\right],
	\end{equation}
	where $u_k= \frac{F(z_k)^\top(x_k - {z_k})}{\|F(z_k)\|^2}$. Let $k :=k +1$, and go to Step 1.
	
	\KwOut{$x_k$.}
\end{algorithm}
%\vspace{-0.2cm}

Algorithm \ref{algo1} is a special case of Algorithm UAF \cite{OL23} that adopts the line search scheme VI. We focus on this scenario since it is representative of DFPM. Several general characters of the framework are analyzed as follows.

{\bf Search direction $d_k$.}
The conditions \eqref{SDC} and \eqref{TRP} for $d_k$ are to guarantee the global convergence. If $F$ is the gradient of a function $f: \mathbb{R}^n \to \mathbb{R}$, then \eqref{SDC} indicates that $d_k$ is a sufficient descent direction for $f$ at $x_k$. Further, the condition \eqref{SDC} implies that the line search procedure \eqref{ak} is well-defined. If $\|d_k\|$ is large during the iteration, the right-hand side of \eqref{ak} will be large, which could lead to more function evaluations and thus increased computational cost. The condition \eqref{TRP} gives $d_k$ a vanishing upper bound, and the method can avoid taking steps that are too long.
The way to obtain $d_k$ satisfying \eqref{SDC} and \eqref{TRP} depends on the particular instance of the framework. For example, the directions in \cite{AK15,MJJYH23,IKRPA22,WA22} all satisfy these conditions. Three specific examples are presented in Section \ref{NE}.

{\bf Line search procedure.} Note that ${\eta_k}(i):=P_{[t_1,t_2]}[\|F(x_k + \gamma \rho ^i d_k)\|]$ in right-hand of \eqref{ak} can be replace by other procedures, for instance ${\eta_k}(i)=\lambda_k+(1-\lambda_k)\|F(x_k + \gamma \rho ^i d_k)\|$, $\lambda_k\in(0,1]$ in \cite{OL18}. Here we only focus on this case since it is a adaptive line search procedure recently proposed by Yin et al. \cite{YJJLW21} and is widely used to compute a stepsize \cite{MJJYH23,WSLZ23}. More specifically, if $t_1=t_2= 1$, then ${\eta_k}(i) = 1$, and thus it reduces to the procedure in \cite{ZZ06}; If $t_1 = 0$ and $t_2$ is large enough, then ${\eta_k}(i)=\|F(x_k + \gamma \rho ^i d_k)\|$, and thus it reduces to the procedure in \cite{LL11}. The projection technique in \eqref{ak} prevents the right-hand side of \eqref{ak} from being too small or too large, which effectively reduces the computational cost of Step 3.

{\bf Projection strategy.} The relaxation factor $\zeta\in (0, 2)$ in \eqref{xk1} serves as a parameter that can enhance the convergence, as stated in \cite{CGH14}. When $\zeta=1$, it corresponds to the original strategy presented in \cite{SS98}. The projection from $x_k$ onto the hyperplane $H_k$ actually provides a descending direction for $x_k$. Although $x_k - \zeta u_k F(z_k)$ is not on $H_k$ for $\zeta \ne 1$, it is still in the direction. ($\zeta$ can be viewed as a stepsize here). As a byproduct of the numerical experiments, we find that taking a suitable relaxation factor $\zeta\in (1,2)$ in the projection step \eqref{xk1} of DFPM can achieve faster convergence.

\subsection{Global Convergence}
We present two simple results to show the global convergence of Algorithm \ref{algo1}. Based on \eqref{SDC} and \eqref{TRP}, the proofs are similar to those of the results in corresponding literature, so we list the results without proof.
\begin{lemma}\cite[Lemma 4]{YJJLW21}\label{lemma3}
	Suppose the sequences $\{x_{k}\}$ and $\{z_{k}\}$ are generated by Algorithm \ref{algo1}. Then the following two claims hold.\\
	(i) For any $x^{*} \in \mathcal{S}$, $\{\|x_k-x^*\|\}$ is convergent.\\
	(ii) $\{x_{k}\}$, $\{d_{k}\}$ and $\{z_{k}\}$ are all bounded, and $\lim_{k \to \infty} \alpha_k  \| d_k \| =0$.
\end{lemma}
\begin{theorem}\cite[Theorem 3.6]{OL23}\label{The1}
	Let sequence $\{x_{k}\}$ be generated by Algorithm \ref{algo1}. Then the sequence $\{x_{k}\}$ converges to a solution of problem (\ref{P}).
\end{theorem}

\section{Anderson Acceleration for DFPM}\label{AAA}
Having seen the convergence for the underlying algorithm, we proceed to show how Anderson acceleration (AA) may translate the improve convergence behavior for DFPM.

\subsection{Anderson Acceleration}
Let $G:\mathbb{R}^n\to \mathbb{R}^n$ be a mapping and consider the problem of finding a fixed-point of $G$:
\[{\rm Find}~x\in\mathbb{R}^n~{\rm such~that}~x=G(x).\]
AA is an efficient acceleration method for fixed-point iteration $x_{k+1} = G(x_k)$. The key idea of AA is to form a new extrapolation point by using the past iterates. To generate a better iterate $x_{k+1}$, it searches for a point $\bar{x}_k$ that has the smallest residual within the subspace spanned by the $m + 1$ most recent iterates. Let $\bar{x}_k=\sum_{j=k-m}^k a^k_j x_j$, $m\leq k$ and $\sum_{j=k-m}^k a^k_j =1$, AA seeks to find a vector of coefficients $a^k=(a^k_{k-m},\dots,a^k_k)^\top$ such that
\[a^k=\arg\min\|G(\bar{x}_k)-\bar{x}_k\|.\]
However, it is hard to find $a^k$ for a general nonlinear mapping $G$. AA uses
\[G(\bar{x}_k)=G\left(\sum_{j=k-m}^k a^k_j x_j\right)\doteq \sum_{j=k-m}^k a^k_j G(x_j),\]
where
\[a^k=\arg\min\left\|\sum_{j=k-m}^k a^k_j G(x_j) -\sum_{j=k-m}^k a^k_j x_j \right\|
=\arg\min\left\|\sum_{j=k-m}^k a^k_j r_j \right\|,\]
with $r_k=G(x_k)-x_k$ to perform an approximation. While $a^k$ is computed, the next iterate of AA is then generated by the following mixing with $b_k\in (0,1]$,
\[
x_{k+1}=(1-b_k)\sum_{j=k-m}^k a^k_j x_j+b_k\sum_{j=k-m}^k a^k_j G(x_{j}).
\]
A formal algorithmic description of AA with the window of length $m_k$ is given by Algorithm \ref{algo2}.

\begin{algorithm}[h]\label{algo2}
	\caption{Anderson acceleration (AA)}
	\SetAlgoLined
	\KwIn{initial point $x_0$, parameters $m>0$, $b_k\in (0,1]$.}
	$x_{1} = G(x_0)$\;
	\For{$k=1,\dots,K-1$}{
		$m_k=\min\{m,k\}$, $r_k=G(x_k)-x_k$\;
		$R_k=(r^k_{k-{m_k}},\dots,r^k_k)^\top$\;
		Solve $\min\limits_{a^k=(a^k_{k-{m_k}},\dots,a^k_k)^\top}\left \|R_k^\top a^k\right\|^2$, subject to $\sum\limits_{j=k-{m_k}}^k a^k_j =1$\;
		$x_{k+1}=(1-b_k)\sum\limits_{j=k-{m_k}}^k a^k_j x_j+b_k\sum\limits_{j=k-{m_k}}^k a^k_j G(x_{j}).$
	}
	\Return{$x_K$.}
\end{algorithm}

In each iteration in Algorithm \ref{algo2}, AA incorporates useful information from previous $m_k$ iterates by an affine combination, where the coefficient $a^k$ is computed as the solution of a minimization problem, rather than expending evaluation directly at current iterate.
One could use any norm in the minimization problem. Using different norms does not affect the convergence. Typically one uses the $\ell_2$ norm, which is what we use here. The reader may refer to \cite{TK15} and references therein for its efficient implementations.
The window size $m$ indicates how many history iterates will be used in the algorithm and its value is typically no larger than 10 in practice. If $m = 0$, AA reduces to the fixed-point iteration. When $b_k$ is a constant independent of $k$, Algorithm \ref{algo2} is referred to as stationary AA. Many works \cite{BC22,MJ20,TK15,WN11} take $b_k\equiv 1$ to simplify the analysis. Here we consider a nonstationary case, and the expression of $b_k$ is given in \eqref{bk}.

\subsection{Acceleration Algorithm}
Based on the convergence result in Section \ref{Alg}, we incorporate AA into DFPM and give the resulting algorithm, named AA-DFPM, in Algorithm \ref{algo3}. Note that a DFPM iteration may not be a fixed-point iteration for $x_k$ since the direction $d_k$ may involve other parameters. However, since AA is a sequence acceleration technique, we expect DFPM to gain a speedup as long as it is convergent.
\begin{algorithm}[h]\label{algo3}
	\caption{AA-DFPM for \eqref{P}}
	\KwIn{initial point $x_0\in\mathcal{C}$, parameters $m,~c,~\gamma,~\sigma,~\epsilon>0$, $0<s_1\leq s_2$, $\rho\in (0,1)$, $0\leq t_1\ll t_2$, $\zeta\in(0,2)$, $b_k\in (0,1]$. Set $k:=0$.}
	
	{\bf Step 1.} Compute $F_k$. If $\|F_k\|<\epsilon $, stop. Otherwise, go to Step 2.
	
	{\bf Step 2.} Compute the search direction $d_k$ satisfying \eqref{SDC} and \eqref{TRP}.
	
	{\bf Step 3.} Choose the stepsize ${\alpha_k}$ satisfying \eqref{ak}, and set $z_k = x_k+\alpha_k d_k$.
	
	{\bf Step 4.} Calculate
	\[v_{k } = P_\mathcal{C}\left[ x_k - \zeta u_k F(z_k)\right],\]
	where $u_k=\frac{F(z_k)^\top(x_k - {z_k})}{\|F(z_k)\|^2}$. If $\|F(v_{k } )\|<\epsilon $, stop. Otherwise, go to Step 5.
	
	{\bf Step 5.} Anderson acceleration for $k\ne0$: set $m_k=\min\{m,k\}$, $r_k=v_k -x_k$. Let $a^k=(a^k_{k-{m_k}},\dots,a^k_k)^\top$, $R_k=(r^k_{k-{m_k}},\dots,r^k_k)^\top$, and solve
	\begin{equation}\label{comak}
		\min\limits_{a^k}\left \|R_k^\top a^k\right\|^2,~{\rm subject~to} \sum\limits_{j=k-{m_k}}^k a^k_j =1,~ a^k_j\geq 0,~j=k-{m_k},\dots,k.
	\end{equation}
	\begin{equation}\label{AA}
		x_{k}^{AA}=(1-b_k)\sum\limits_{j=k-{m_k}}^k a^k_j x_j+b_k\sum\limits_{j=k-{m_k}}^k a^k_j v_{j}.
	\end{equation}
	If 
	\begin{equation}\label{sgc}
		\left\|\sum\limits_{j=k-{m_k}}^{k}a^k_j x_j-v_k\right\|\leq ck^{-(1+\epsilon)},
	\end{equation}
	then $x_{k+1}=x_{k}^{AA}$, else $x_{k+1}=v_{k}$.
	Let $k:=k +1$, and go to Step 1.
	
	\KwOut{$x_k$.}
\end{algorithm}

Some implementation techniques in the algorithm bear further commenting. We thus introduce and discuss the following four aspects.

\noindent{\bf Feasibility of accelerated iterate.} As illustrated in Algorithm \ref{algo2}, AA computes the accelerated iterate via an affine combination of previous iterates. The accelerated point may violate the constraint unless its feasible set is affine. Considering that the feasible set $\mathcal{C}$ in problem \eqref{P} is closed and convex and the previous iterates generated by Algorithm \ref{algo1} are all in $\mathcal{C}$, we set $a^k_j\geq 0,~j=k-{m_k},\dots,k$, in \eqref{comak} to obtain a reliable accelerated iterate. This means that the accelerated iterate here is computed by a convex combination of previous iterates. Version to this technique is called EDIIS in the chemistry community \cite{KSC02}.

\noindent{\bf Computation of coefficient $a^k$.} The residual matrix $R_k$ in the least squares problem of AA can be rank-deficient; then ill-conditioning may occur in computing $a^k$. Here a Tikhonov regularization \cite{SDB16} $\lambda\|a^k\|^2,~\lambda>0$, be added to the problem to obtain a reliable $a^k$. In this case, combining with the above implementation, the least squares problem is a standard quadratic programming problem that can be solved by the MATLAB command ``quadprog''.

\noindent{\bf Guarantee of convergence.} As mentioned earlier, since AA is known to only converge locally, some globalization mechanisms are required to use it in practice, such as adaptive regularization \cite{OTMD23}, restart checking \cite{HV19} and safeguarding step \cite{OPYZD20}. Following \cite{ZOB20}, we introduce a safeguard checking \eqref{sgc} to ensure the global convergence, thus
\begin{equation}\label{axv}
	\sum_{k=1}^{\infty} \left\|\sum_{j=k-{m_k}}^{k}a^k_j x_j-v_k\right\|<\infty.
\end{equation}

\noindent{\bf Calculation of $b_k$.} Define the following averages with the solution $a^k$ to the least square problem in Step 5 of Algorithm \ref{algo3},
\[x_k^{a}=\sum_{j=k-{m_k}}^k a^k_j x_j,~v_k^{a}=\sum_{j=k-{m_k}}^k a^k_j v_{j}.\]
Then \eqref{AA} becomes
\begin{equation}\label{AAI}
	x_{k}^{AA}=(1-b_k)x_k^{a}+b_k v_k^{a}=x_k^{a}+b_k(v_k^{a}-x_k^{a}).
\end{equation}
The relaxation parameter $b_k$ is generally determined heuristically. Many discussions choose $b_k\equiv 1$, thereby simplifying the expression to facilitate theoretical analysis. Little attention has been paid to nonstationary case. As Anderson wished in his comment \cite{Anderson19}, we design a dynamic factor
\begin{equation}\label{bk}
	b_k=\min\left\{b,\frac{1}{k^{(1+\epsilon)}\| v_k^{a}-x_k^{a}\|}\right\},~b\in(0,1).
\end{equation}
The adaptive idea is derived from the inertial-based algorithms \cite{CMY14,MJJYH23}. Then for all $k$, we have $b_{k}\|v_k^{a}-x_k^{a}\| \leq k^{-(1+\epsilon)}$, which implies that
\begin{equation}\label{int}
	\sum_{k=1}^{\infty} b_{k}\|v_k^{a}-x_k^{a}\|<\infty.
\end{equation}
\begin{remark}
	A major difference in the acceleration strategies between the two schemes: our method is an interpolation procedure that uses a convex combination of iterates, whereas the original AA is actually an extrapolation procedure that uses an affine combination of iterates.
\end{remark}

\section{Convergence Analysis}\label{Cov}
We first present the following lemma to help us complete the proof. 
\begin{lemma}\cite{MJJYH23}\label{lemma4}
	Let $\{ \alpha_k  \} $ and $\{ \beta_k  \} $ be two sequences of nonnegative real numbers satisfying $\alpha_{k+1} \leq \alpha_{k}+\beta_{k}$ and $\sum_{k=1}^{\infty} \beta_{k}<+\infty$. Then the sequence $\{ \alpha_k  \} $ is convergent as $k\to \infty$.
\end{lemma}

This lemma is derived from \cite[Lemma 9]{Polyak1987}, which is a result on random variables. The proof of Lemma \ref{lemma4} has been proven in \cite[Lemma 2]{MJJYH23}, so its proof is omitted here.

We can now get the convergence results for AA-DFPM. 
\begin{lemma}\label{lemma5}
	Let sequence $\{x_{k}\}$ be generated by Algorithm \ref{algo3}. Then for any $x^{*} \in \mathcal{S}$, (i) the sequence $\{\|x_{k}-x^*\|\}$ is convergent; (ii) $\lim_{k \to \infty} \alpha_k  \| d_k \| =0$.
\end{lemma}
{\it Proof}
Depending on whether the sequence processes AA or not, we partition the iteration counts into two subsets accordingly, with $K_{AA} = \{ k_0, k_1,\dots\}$ being those iterations passing \eqref{sgc} and $K_{O} = \{ l_0, l_1,\dots\}$ being the rest.

Consider $x^* \in \mathcal{S}$ a solution of \eqref{P}. In the following derivation, we assume that both $K_{AA}$ and $K_{O}$ are infinite. The cases when either of them is finite are even simpler as one can completely ignore the finite index set.

(i) For $l_i\in K_{O}\ (i>0)$, by inequality (19) in \cite{YJJLW21}, we have that
\begin{equation}\label{xw*}
	\|x_{l_i+1}-x^{*}\|^{2}\leq \|x_{l_i}-x^{*}\|^{2}-\zeta (2-\zeta)\frac{\sigma^{2}t_1^2 \alpha_{l_i}^4\| d_{l_i} \|^4}{\|F(z_{l_i})\|^{2}}\leq \|x_{l_i}-x^{*}\|^{2}.
\end{equation}

For $k_i\in K_{AA}\ (i>0)$, from \eqref{AAI}, we have
\[\begin{aligned}
	\|x_{k_i+1}-x^{*}\| &= \|x_{k_i}^{a}+b_{k_i}(v_{k_i}^a-x_{k_i}^{a})-x^{*}\|
	 \leq \|x_{k_i}^{a}-x^{*}\|+b_{k_i}\|v_{k_i}^{a}-x_{k_i}^{a}\|\\
	&\leq \| v_{k_i}-x^{*}\|+\left\|\sum_{j={k_i}-{m_{k_i}}}^{k_i} a^{k_i}_j x_j-v_{k_i}\right\|+b_{k_i}\|v_{k_i}^{a}-x_{k_i}^{a}\|.
\end{aligned}\]
Similar to the proof of \eqref{xw*}, we can get
\begin{equation}\label{vx*}
	\|v_{k_i}-x^{*}\|^{2}\leq \|x_{k_i}-x^{*}\|^{2}-\zeta (2-\zeta)\frac{\sigma^{2}t_1^2 \alpha_{k_i}^4\| d_{k_i} \|^4}{\|F(z_{k_i})\|^{2}}\leq \|x_{k_i}-x^{*}\|^{2}.
\end{equation}
Hence
\begin{align}
	\|x_{k_i+1}-x^{*}\|&\leq \|x_{k_i}-x^{*}\|+\left\|\sum\limits_{j={k_i}-{m_{k_i}}}^k a^{k_i}_j x_j-v_{k_i}\right\|+b_{k_i}\|v_{k_i}^{a}-x_{k_i}^{a}\|\nonumber\\
	& \leq \|x_{k_i}-x^{*}\|+\beta_i,\label{xki}
\end{align}
with $\beta_i=(1+c) i^{-(1+\epsilon)}$. 

By telescoping \eqref{xw*} and \eqref{xki}, we obtain that
\[\|x_{k+1}-x^{*}\|\leq \|x_k-x^{*}\|+\beta_k,\]
with $\beta_k\geq 0$ and $\sum_{k=0}^\infty \beta_k<\infty$.
Using Lemma \ref{lemma4} with $\alpha_{k}=\|x_{k}-x^{*}\|$, the sequence $\{\|x_{k}-x^*\|\}$ is convergent. 

(ii) The above result implies that $\{x_k\}$ is bounded. This, together with the continuity of $F$ and \eqref{TRP}, shows that $\{d_k\}$ is bounded, further implies $\{z_k\}$ is bounded, as well as $\{F(z_k)\}$.
Suppose $\|F(z_k)\|\leq N$ and $\|x_{k}-x^*\|\leq M$. Summing \eqref{xw*}, we have
\[
\begin{aligned}
	\zeta(2-\zeta)\frac{\sigma^{2}t_1^2}{N^{2}}\sum\limits_{i=0}^{\infty} (\alpha_{l_i}\| d_{l_i} \|)^{4} &\leq\sum\limits_{i=0}^{\infty}(\|{x}_{{l_i}}-x^{*}\|^{2}
	-\|x_{{l_i}+1}-x^{*}\|^{2}) \\
	&=\|{x}_{0}-x^{*}\|^{2}-\lim_{i \to \infty} \|{x}_{{l_i}+1}-x^{*}\|^{2}<\infty.
\end{aligned}
\]
Hence $ \lim_{i \to \infty} \alpha_{l_i}  \| d_{l_i} \| =0$.

On the other hand,
\begin{align}
	\|x_{k_i}^{a}-x^{*}\|&=\left\|\sum_{j={k_i}-{m_{k_i}}}^{k_i} a^{k_i}_j (x_j-x^{*})\right\|\leq \sum_{j={k_i}-{m_{k_i}}}^{k_i} |a^{k_i}_j| \|x_j-x^{*}\|\nonumber\\
	&\leq (m_k+1) M \leq (m+1) M.\label{xaM}
\end{align}
By \eqref{sgc} and $\|v_{k_i}-x^{*}\|\leq \|x_{k_i}-x^{*}\|\leq M$, we have
\begin{align}
	\|x_{k_i}^{a}-x^{*}\|^2&\leq (\| v_{k_i}-x^{*}\|+\|x_{k_i}^{a}-v_k\|)^2\leq (\| v_{k_i}-x^{*}\|+c)^2\nonumber\\
	&=\| v_{k_i}-x^{*}\|^2+2c\|v_{k_i}-x^{*}\|+c^2\nonumber\\
	&\leq\| v_{k_i}-x^{*}\|^2+2cM+c^2.\label{xa2M}
\end{align}
Therefore, we get
\[
\begin{aligned}
	\|x_{{k_i}+1}-x^{*}\|^2 \leq & \|x_{k_i}^{a}-x^{*}\|^2+(b_{k_i}\|v_{k_i}^{a}-x_{k_i}^{a}\|)^2+2b_{k_i}\|v_{k_i}^a-x_{k_i}^{a}\| \|x_{k_i}^{a}-x^{*}\|\\
	\overset{\eqref{xaM}}{\leq} & \|x_{k_i}^{a}-x^{*}\|^2+k_i^{-(2+2\epsilon)}+2(m+1)Mk_i^{-(1+\epsilon)}\\
	\overset{\eqref{xa2M}}{\leq} &\| v_{k_i}-x^{*}\|^2+2cM+c^2+k_i^{-(2+2\epsilon)}+2(m+1) Mk^{-(1+\epsilon)}\\
	\overset{\eqref{vx*}}{\leq} &\|x_{k_i}-x^{*}\|^{2}-\zeta(2-\zeta)\frac{\sigma^{2}t_1^2 \alpha_{k_i}^4\| d_{k_i} \|^4}{\|F(z_{k_i})\|^{2}}+2cM+c^2\\
	&+k_i^{-(2+2\epsilon)}+2(m+1) M k^{-(1+\epsilon)}.
\end{aligned}
\]
Adding above inequality, in view of the boundedness of $\{F(z_k)\}$ and \eqref{int}, it follows that
\[
\begin{aligned}
	&\zeta (2-\zeta )\frac{\sigma^{2}t_1^2}{N^{2}}\sum_{i=0}^{\infty}(\alpha_{k_i}\| d_{k_i} \|)^{4}\\
	\leq&\sum_{i=0}^{\infty}\left[\|{x}_{k_i}-x^{*}\|^{2}
	-\|x_{{k_i}+1}-x^{*}\|^{2}+2cM+c^2+k_i^{-(2+2\epsilon)}+2(m+1) M k^{-(1+\epsilon)} \right]\\
	=&\|{x}_{0}-x^{*}\|^{2}-\lim_{i \to \infty} \|{x}_{{k_i}+1}-x^{*}\|^{2}+2cM+c^2+\sum_{i=0}^{\infty}k_i^{-(2+2\epsilon)}\\
	&+2(m+1) M \sum_{i=0}^{\infty}k^{-(1+\epsilon)} <\infty.
\end{aligned}
\]
This implies $\lim_{i \to \infty} \alpha_{k_i}  \| d_{k_i} \|=0$. 
Together with $ \lim_{i \to \infty} \alpha_{l_i}  \| d_{l_i} \| =0$, we have $\lim_{k \to \infty} \alpha_k  \| d_k \|=0$.
\qed

\begin{remark}
	The monotonicity of $F$ is a common assumption \cite{LL11,MJJYH23,WA22,YJJLW21} for DFPM to construct the hyperplane (see \eqref{mon}) whose projection provides a descending direction for $x_k$. It is also essential to obtain the descent of the sequence $\{\|x_{k}-x^*\|\}$ (i.e. \eqref{xw*}) in its convergence analysis. Through fixed-point mappings or normal mappings \cite{ZL01}, a number of monotone variational inequality problems can be converted into monotone systems. Some sufficient conditions for their monotonicity have been discussed in \cite{ZL01}. In addition, some works have explored the new DFPM whose $F$ is pseudo-monotonicity \cite{JYTH22,LLXWT22}.
\end{remark}

Based on Lemma \ref{lemma5}, we prove the global convergence for AA-DFPM.
\begin{theorem}\label{The2}
	The sequence $\{x_{k}\}$ generated by Algorithm \ref{algo3} converges to a solution of problem (\ref{P}).
\end{theorem}
{\it Proof}
	Assume that $\liminf\limits_{k \to \infty}\|F_k\|>0$, there exists a constant $\varepsilon>0$ such that
	\begin{equation}\label{1.21}
		\|F_k\| \geq \varepsilon, \forall~k \geq 0.
	\end{equation}
	Further, from \eqref{SDC} and Cauchy-Schwarz inequality, we have
	\[
	\|d_{k}\| \geq s_1\|F_k\| \geq s_1\varepsilon> 0,~\forall~k \geq 0.
	\]
	This together with Lemma \ref{lemma5} (ii) implies
	\begin{equation}\label{1.22}
		\lim _{k \to \infty} \alpha_{k}=0.
	\end{equation}
	In view of the boundedness of $\{x_k\}$ and $\{d_k\}$, there exist two subsequences $\{x_{k_j}\}$ and $\{d_{k_j}\}$ such that
	\[
	\lim\limits _{j \to \infty} x_{k_j}=\hat{x},~\lim\limits _{j \to \infty} d_{k_j}=\hat{d}.
	\]
	Again, it follows from \eqref{SDC} that
	\[
	-F_{k_j}^\top d_{k_j} \geq s_1 \|F_{k_j}\|^{2},~\forall~j.
	\]
	Letting $j \to \infty$ in the inequality above, and by the continuity of $F$ and \eqref{1.21}, we get
	\begin{equation}\label{1.23}
		-F(\hat{x})^\top \hat{d} \geq s_1 \|F(\hat{x})\|^{2} > s_1 \varepsilon^{2}>0.
	\end{equation}
	Similarly, it follows from \eqref{ak} that
	\[
	-F(x_{k_j}+\rho^{-1} \alpha_{k_j} d_{k_j})^\top d_{k_j}<\sigma \rho^{-1} \alpha_{k_j} P_{[t_1,t_2]}[\|F(x_{k_j}+\rho^{-1} \alpha_{k_j} d_{k_j})\|] \|d_{k_j}\|^2,~\forall~j.
	\]
	Letting $j \to \infty$ in the inequality above, taking into account \eqref{1.22} and the continuity of $F$, we conclude that $-F(\hat{x})^\top \hat{d} \leq 0$, which contradicts \eqref{1.23}. Thus,
	\begin{equation}\label{1.20}
		\liminf\limits_{k \to \infty}\|F_k\|=0.
	\end{equation}
	By the boundedness of $\{x_{k}\}$ and the continuity of $F$ as well as \eqref{1.20}, the sequence $\{x_{k}\}$ has an accumulation point $x^*$ such that $F(x^*)=0$. By $x_k \in \mathcal{C}$ and the closeness of $\mathcal{C}$, we have $x^* \in \mathcal{C}$, further $x^* \in \mathcal{S}$. Combining with the convergence of $\{\|x_{k}-x^*\|\}$ (Lemma \ref{lemma5} (i)), one knows that the whole sequence $\{x_{k}\}$ converges to $x^*\in \mathcal{S}$.
\qed

By Theorem \ref{The2}, we can assume that $x_k \to x^*\in \mathcal{S}$ as $k \to \infty$. Under mild assumptions below, we further illustrate the linear convergence rate of AA-DFPM.

\begin{assumption}\label{ass1}
	The mapping $F$ is Lipschitz continuous on $\mathbb{R}^n$, i.e., there exists a positive constant $L$ such that
	\begin{equation}\label{LC}
		\| F(x)-F(y) \|\leq L\| x-y \| ,~\forall~x,y \in \mathbb{R}^n.
	\end{equation}
\end{assumption}

The Lipschitz continuity assumption on $F$ helps us to provide a uniform lower bound of the stepsize $\alpha_{k}$. Based on \eqref{SDC}, \eqref{TRP} and \eqref{LC}, the proof is similar to that of Lemma 3.4 in \cite{WSLZ23}, so we omit it here.
\begin{lemma}\label{lemma1}
	Suppose that Assumption \ref{ass1} holds. Then the stepsize $\alpha_k$ yielded by (\ref{ak}) satisfies
	\begin{equation}\label{1.3}
		\alpha_{k}\geq\alpha :={\rm min}\left \{ \gamma ,~\frac{\rho s_1}{(L+\sigma t_2)s_2^2}\right \} > 0.
	\end{equation}
\end{lemma}

\begin{assumption}\label{ass2}
	For the limit $x^*\in \mathcal{S}$ of $\{x_k\}$, there exist two positive constants $\ell $ and $\varepsilon $ such that,
	\begin{equation}\label{1.13}
		\ell~{\rm dist}(x_k,\mathcal{S})\leq \| F_k \|,~\forall~x_k \in B(x^*,\varepsilon),~k=1,2,\dots,
	\end{equation}
	where the neighborhood $B(x^*,\varepsilon)=\{ x_k \in \mathbb{R}^n:\ \| x_k-x^*  \| <\varepsilon \}$.
\end{assumption}

The local error bound Assumption \ref{ass2} is usually used to prove the convergence rate of DFPM in solving \eqref{P} (see \cite{LF19,MJJYH23,OL18} for instance). It holds whenever constrained set $\mathcal{C}$ is polyhedral and either function $F$ is affine or $F$ is strongly monotone and Lipschitz continuous on $\mathcal{C}$ (see Theorem 2.2 in \cite{T95}). Now we estimate the asymptotic rate of convergence of the iteration, for sufficiently large $k$. The sequence $\{ x_k \}$ in the proof of the following theorems refers to the acceleration iteration. The convergence rate of the original iteration is identical to Theorem 4.5 in \cite{OL23}.

\begin{theorem}\label{thm3}
	Suppose that Assumptions \ref{ass1} and \ref{ass2} hold, and the sequence $\{ x_k \}$ is generated by Algorithm \ref{algo3}. Then $ \{ {\rm dist}(x_k,\mathcal{S}) \}$ satisfies
	\[
	\frac{{\rm dist}(x_{k+1}, \mathcal{S})}{{\rm dist}(x_{k}, \mathcal{S})} \leq \sqrt{\varphi}+(c+1)\frac{k^{-(1+\epsilon)}}{{\rm dist}(x_{k}, \mathcal{S})},
	\]
	where $\varphi  = 1-\zeta (2-\zeta )\left(\frac{\sigma t_1 \alpha^2 s_1^2\ell^2}{\varrho}\right)^2$ and $\varrho=\max\{L(\gamma L s_2+1),\sqrt{\zeta(2-\zeta)}\sigma t_1 \alpha^2 s_1^2\ell^2\}$.
\end{theorem}
{\it Proof}
Let $h_k \in \mathcal S$ be the closest solution to $x_k$, i.e., $\left \| x_k-h_k \right \|={\rm dist}(x_k,\mathcal S)$. Recall \eqref{vx*} that
\begin{align}
	\|v_{k}-h_{k}\|^{2}&\leq\|x_{k}-h_{k}\|^{2}-\zeta(2-\zeta )\frac{\sigma^{2}t_1^2 \alpha_k^4\| d_k \|^4}{\|F(z_{k})\|^{2}}\nonumber\\
	&={\rm dist}^{2}(x_{k}, \mathcal{S})-\zeta(2-\zeta)\frac{\sigma^{2}t_1^2 \alpha_k^4\| d_k \|^4}{\|F(z_{k})\|^{2}}.\label{1.15}
\end{align}
From \eqref{LC}, \eqref{TRP} and $0< \alpha_{k} \leq \gamma$, it follows from that
\begin{align}
	\|F(z_{k})\| &=\|F(z_{k})-F(h_{k})\| \overset{\eqref{LC}}{\leq} L\|z_{k}-h_{k}\|\leq L(\|x_{k}-z_{k}\|+\|x_{k}-h_{k}\|) \nonumber\\
	&=L(\alpha_{k}\|d_{k}\|+\|x_{k}-h_{k}\|) \leq L(\gamma\|d_{k}\|+\|x_{k}-h_{k}\|) \nonumber\\
	& \overset{\eqref{TRP}}{\leq} L(\gamma s_2\|F_{k}\|+\|x_{k}-h_{k}\|) =L(\gamma s_2\|F_{k}-F(h_k)\|+\|x_{k}-h_{k}\|) \nonumber\\
	& \overset{\eqref{LC}}{\leq} L(\gamma L s_2+1)\|x_{k}-h_{k}\|=L(\gamma L s_2+1){\rm dist}(x_k,\mathcal S)\nonumber\\
	&\leq \varrho {\rm dist}(x_k,\mathcal S).\label{1.16}
\end{align}
Again, from \eqref{SDC}, \eqref{1.3} and \eqref{1.13}, we have
\begin{equation}\label{1.17}
	\alpha_k^4\| d_k \|^4\geq \alpha^4 s_1^4\|F_{k}\|^4\overset{\eqref{1.13}}{\geq} \alpha^4 s_1^4\ell^4{\rm dist}^4(x_{k}, \mathcal{S}).
\end{equation}
Combining with (\ref{1.15})-(\ref{1.17}), we obtain
\[
	\|v_{k}-h_{k}\|^{2}\leq \varphi  {\rm dist}^{2}(x_{k}, \mathcal{S}).
\]
This, together with $\|x_k^{a}-v_k\|\leq ck^{-(1+\epsilon)}$ and $b_k\|v_k^{a}-x_k^{a}\| \leq k^{-(1+\epsilon)}$, shows that
\begin{align}
	{\rm dist}(x_{k+1}, \mathcal{S}) & \leq \|x_{k+1}-h_k\| \nonumber\\
	&=\|v_k- h_k+x_k^{a}-v_k+b_k(v_k^{a}-x_k^{a})\| \nonumber\\
	&\leq \|v_k- h_k\|+\|x_k^{a}-v_k\|+b_k\|v_k^{a}-x_k^{a}\| \nonumber\\
	&\leq \sqrt{\varphi}{\rm dist}(x_{k}, \mathcal{S})+(c+1)k^{-(1+\epsilon)}. \label{dk1dk}
\end{align}
Hence
\[
\frac{{\rm dist}(x_{k+1}, \mathcal{S})}{{\rm dist}(x_{k}, \mathcal{S})} \leq \sqrt{\varphi}+(c+1)\frac{k^{-(1+\epsilon)}}{{\rm dist}(x_{k}, \mathcal{S})}.
\]
The proof is completed.
\qed

Let $a:=\limsup\limits_{k \to \infty}\frac{k^{-(1+\epsilon)}}{{{\rm dist}(x_k,\mathcal S)}}$. From Theorem \ref{thm3}, the existence of $a$ is essential to further obtain the convergence rate results. Different $a$ correspond to different convergence rates of the sequence $ \{{\rm dist}(x_k,\mathcal{S}) \}$. Its value provides insight into the following asymptotic behavior.
\begin{corollary}
	Suppose that Assumptions \ref{ass1} and \ref{ass2} hold, and the sequence $\{ x_k \}$ is generated by Algorithm \ref{algo3}. Then the three following claims hold.
	
	\noindent(i)~If $0<a<+\infty $, then ${{\rm dist}(x_{k},\mathcal S)}=O (k^{-(1+\epsilon)})$;
	
	\noindent(ii)~If $a=+\infty $, then ${{\rm dist}(x_{k},\mathcal S)}=o (k^{-(1+\epsilon)})$;
	
	\noindent(iii)~If $a=0$, then the sequence $ \{ {\rm dist}(x_k,\mathcal{S}) \}$ converges Q-linearly to 0, i.e.,
	\[
	\limsup\limits_{k \to \infty}\frac{{\rm dist}(x_{k+1},\mathcal S)}{{\rm dist}(x_k,\mathcal S)} <1.
	\]
\end{corollary}

This result is consistent with the convergence rate of the inertial-type DFPM \cite{MJJYH23,WSLZ23}. We further investigate the convergence rate of sequence $\{x_k\}$ if the mapping $F$ is strongly monotone with modulus $\mu > 0$, i.e.,
\[
(F(x)-F(y))^{\top}(x-y)\geq \mu\|x-y\|^2,~\forall~x,~y\in \mathbb{R}^n.
\]
\begin{theorem} \label{thm4}
Suppose that Assumptions \ref{ass1} and \ref{ass2} hold, and the sequence $\{ x_k \}$ is generated by Algorithm \ref{algo3}. If the mapping $F$ is strongly monotone, then $ \{ \|x_k-x^*\|\}$ satisfies
	\[\frac{\|x_{k+1}-x^*\|}{\|x_k-x^*\|} \leq \sqrt{\psi}+(c+1)\frac{k^{-(1+\epsilon)}}{\|x_k-x^*\|},\]
	where $\psi=1-\zeta (2-\zeta) \left(\frac{\sigma t_1 \alpha^2 s_1^2\mu^2}{\xi}\right)^2$ and $\xi=\max\{L(\gamma L s_2+1),\sqrt{\zeta (2-\zeta)}\sigma t_1 \alpha^2 s_1^2\mu^2\}$.
\end{theorem}
{\it proof}
	By the Cauchy-Schwarz inequality and the strong monotonicity of $F$, it has
	\[\|F_k\|=\|F_k-F(x^*)\|\geq \mu\|x_k-x^*\|.\]
	Together with \eqref{SDC} and \eqref{1.3}, we have
	\begin{equation}\label{xcr1}
		\alpha_k^4\| d_k \|^4 \geq \alpha^4 s_1^4\|F_{k}\|^4\geq \alpha^4 s_1^4\mu^4\|x_k-x^*\|^4.
	\end{equation}
	Similar to the proof of \eqref{1.16}, it follows that
	\begin{equation}\label{xcr2}
		\|F(z_{k})\| \leq L(\gamma L s_2+1)\|x_k-x^*\|\leq \xi \|x_k-x^*\|.
	\end{equation}
	Combining \eqref{vx*} with \eqref{xcr1} and \eqref{xcr2} implies
	\[
		\|v_{k}-x^*\|^{2}\leq \psi\|x_k-x^*\|^2.
	\]
	Also similar to the proof of \eqref{dk1dk}, we obtain
	\[
	\begin{aligned}
		\|x_{k+1}-x^*\|&=\|v_k- x^*+x_k^{a}-v_k+b_k(v_k^{a}-x_k^{a})\|\\
		&\leq \|v_k- x^*\|+\|x_k^{a}-v_k\|+b_k\|v_k^{a}-x_k^{a}\| \\
		&\leq \sqrt{\psi} \|x_k-x^*\|+(c+1)k^{-(1+\epsilon)}.
	\end{aligned}
	\]
	Thus
	\[
	\frac{\|x_{k+1}-x^*\|}{\|x_k-x^*\|} \leq \sqrt{\psi}+(c+1)\frac{k^{-(1+\epsilon)}}{\|x_k-x^*\|}.
	\]
	The proof is completed.
\qed

Let $A:=\limsup\limits_{k \to \infty}\frac{k^{-(1+\epsilon)}}{\|x_k-x^*\|}$. We can also get the asymptotic convergence rate of sequence $\{ x_k \}$ from Theorem \ref{thm4}.
\begin{corollary}
Suppose that Assumptions \ref{ass1} and \ref{ass2} hold, and the sequence $\{ x_k \}$ is generated by Algorithm \ref{algo3}. Then the following statements hold.
	
	\noindent(i)~If $0<A<+\infty $, then $\|x_k-x^*\|=O (k^{-(1+\epsilon)})$;
	
	\noindent(ii)~If $A=+\infty $, then $\|x_k-x^*\|=o (k^{-(1+\epsilon)})$;
	
	\noindent(iii)~If $A=0$, then the sequence $ \{x_k\}$ linearly converges to $x^*\in \mathcal{S}$, i.e.,
	\[\limsup\limits_{k \to \infty}\frac{\|x_{k+1}-x^*\|}{\|x_k-x^*\|} <1.\]
\end{corollary}
\section{Numerical Experiments and Applications}\label{NE}
In this section, we demonstrate the effectiveness of AA-DFPM through experiments on constrained nonlinear equations as well as a real-world problem of machine learning. All tests are conducted in MATLAB R2016b on a 64-bit Lenovo laptop with Intel(R) Core(TM) i7-6700HQ CPU (2.60 GHz), 16.00 GB RAM and Windows 10 OS. Throughout the numerical experiments, three search directions are chosen as follows:

\noindent 1) {\it Spectral conjugate gradient projection (SCGP) method} \cite{WSLZ23}
\begin{equation}\label{SCGP}
	d_k=\left\{\begin{array}{ll}
		-F_k,&k=0,\\
		-\theta_{k}F_k+\beta_k d_{k-1},&k\geq1~{\rm and}~\theta_{k}\in [\vartheta_1,\vartheta_2],\\
		-F_k+\zeta \frac{\|F_k\|}{\|d_{k-1}\|}d_{k-1},&k\geq1~{\rm and}~\theta_{k}\notin [\vartheta_1,\vartheta_2],
	\end{array}\right.
\end{equation}
where
\[\begin{gathered}
	\beta_k=\max \left\{\frac{F_k^{\top} \eta_{k-1}}{d_{k-1}^{\top} v_{k-1}}-\frac{\|\eta_{k-1}\|^2 F_k^{\top} d_{k-1}}{(d_{k-1}^{\top} v_{k-1})^2}, \chi \frac{F_k^{\top} d_{k-1}}{\|d_{k-1}\|^2}\right\},\\
	\eta_{k-1}=y_{k-1}+\tau_k F_k,~\tau_k=\tau \frac{\|y_{k-1}\|}{\|F_k\|}+\min \left\{0, \frac{-F_k^{\top} y_{k-1}}{\|F_k\|^2}\right\},\\
	v_{k-1}=y_{k-1}+\lambda_k d_{k-1},~\lambda_k=\frac{\|y_{k-1}\|}{\|d_{k-1}\|}+\max \left\{0, \frac{-d_{k-1}^{\top} y_{k-1}}{\|d_{k-1}\|^2}\right\},\\
	\theta_k=\frac{s_{k-1}^{\top} F_k+\beta_k y_{k-1}^{\top} d_{k-1}}{F_k^{\top} y_{k-1}},~s_{k - 1}= x_k - x_{k-1},~y_{k - 1}= F_k - F_{k-1},\\
\end{gathered}\]
and $\chi\in(0,\frac{1}{4})$, $\zeta\in[0,1)$, $\tau>0$, $\frac{1}{4}<\vartheta_1<\vartheta_2$. From Lemma 3.1 in \cite{WSLZ23}, direction \eqref{SCGP} satisfies conditions \eqref{SDC} and \eqref{TRP}.

\noindent 2) {\it Hybrid three-term conjugate gradient projection (HTTCGP) method} \cite{YJJLW21}
\begin{equation}\label{HTTCGP}
	d_k=\left\{\begin{array}{ll}-F_k,~&k=0,\\-F_k+\beta_{k} d_{k-1}+\tilde{\upsilon}_k y_{k-1},~&k\geq1, \end{array}\right.
\end{equation}
where
\[\begin{gathered}
	\beta_{k}=\frac{F_{k}^\top y_{k-1}}{\tau_{k}}-\frac{\|y_{k-1}\|^{2} F_{k}^\top d_{k-1}}{\tau_{k}^{2}},~\tilde{\upsilon}_k =\delta_{k} \frac{F_{k}\top d_{k-1}}{\tau_{k}}, \\
	\tau_{k} =\max \{\mu\|d_{k-1}\| \|y_{k-1}\|,~d_{k-1}^\top y_{k-1},~\|F_{k-1}\|^{2}\},
\end{gathered}\]
with parameters $\mu >0$ and $0\leq \delta_k\leq \delta <1$. From Lemma 2 in \cite{YJJLW21}, direction \eqref{HTTCGP} satisfies conditions \eqref{SDC} and \eqref{TRP}.

\noindent 3) {\it Modified spectral three-term conjugate gradient method} \cite{AF23} (Considering that the direction was originally designed for a conjugate gradient method for solving unconstrained problems, here we have adapted it slightly to accommodate DFPM, named MSTTCGP.)
\begin{equation}\label{MSTTCG}
	d_k=\left\{\begin{array}{ll}
		-F_k,&k=0,\\
		-\theta_{k}F_k+\beta_k d_{k-1}- \tilde{\upsilon}_k y_{k-1},&k\geq1~{\rm and}~\theta_{k}\in [\vartheta_1,\vartheta_2],\\
		-F_k+\beta_k d_{k-1}- \tilde{\upsilon}_k y_{k-1},&k\geq1~{\rm and}~\theta_{k}\notin [\vartheta_1,\vartheta_2],
	\end{array}\right.
\end{equation}
where
\[\begin{gathered}
	\beta _k = \frac{F_k^\top y _{k - 1}}{\tau_k},~\tilde{\upsilon}_k=\frac{F_k^\top d _{k - 1}}{\tau_k},~\theta_k=\frac{s_{k-1}^{\top} F_k+\beta_k y_{k-1}^{\top} d_{k-1}-\tilde{\upsilon}_k\|y_{k-1}\|^2}{F_k^{\top} y_{k-1}},\\
	\tau_k = \max \{ \mu \|d_{k - 1}\| \|y_{k - 1}\|,~d_{k - 1}^\top y_{k - 1},~\| F_{k-1} \|^2 \},
\end{gathered}\]
in which $\vartheta_2>\vartheta_1>0$ and $\mu>0$. From Corollary 3.1 in \cite{AF23}, it follows that \eqref{MSTTCG} satisfies condition \eqref{SDC}. We prove that \eqref{MSTTCG} satisfies condition \eqref{TRP}. To proceed, by the definitions of parameters $\theta_{k}$, $\beta_k$ and $\tilde{\upsilon}_k$, we get
\[\begin{aligned}
	\|d_k\| &=\|\theta_k F_k+\beta_k d_{k-1}- \tilde{\upsilon}_k y_{k-1}\| \\
	&\leq\theta_k \|F_k\|+|\beta_k| \|d_{k-1}\|+ |\tilde{\upsilon}_k| \|y_{k-1}\| \\
	&=\theta_k \|F_k\|+\frac{|F_k^\top y _{k - 1}|\|d_{k-1}\|}{\tau_k}+\frac{|F_k^\top d_{k-1}| \|y _{k - 1}\|}{\tau_k}\\
	&\leq \theta_k \|F_k\|+\frac{\|F_k\| \|y _{k - 1}\|\|d_{k-1}\|}{\mu \|d_{k - 1}\| \|y_{k - 1}\|}+\frac{\|F_k\| \|d_{k-1}\| \|y_{k - 1}\|}{\mu \|d_{k - 1}\| \|y_{k - 1}\|}\\
	&=\left(\vartheta_2+\frac{2}{\mu}\right)\|F_k\|,
\end{aligned}\]
for $k\geq1$ and $\theta_{k}\in [\vartheta_1,\vartheta_2]$, and
\[
\|d_k\| =\|\theta_k F_k+\beta_k d_{k-1}- \tilde{\upsilon}_k w_{k-1}\| \leq\left(1+\frac{2}{\mu}\right)\|F_k\|,
\]
for $k\geq1$ and $\theta_{k}\notin [\vartheta_1,\vartheta_2]$. Thus $\|d_k\| \leq s_2 \|F_k\|,~s_2:=\max\{1,\vartheta_2\}+\frac{2}{\mu}$.

All related parameters of SCGP, HTTCGP and MSTTCGP are the same as their originals. In addition, we set the line search and the projection parameters $\sigma=0.01$, $\gamma=1$, $\rho=0.6$, $\xi=1.7$, $t_1=0.001$, $t_2=0.4$ for MSTTCGP. We use AA-SCGP, AA-HTTCGP and AA-MSTTCGP to denote their Anderson acceleration variant with the AA parameters $c=10$, $b=0.1$ and $\lambda=10^{-10}$. We test the effect of $m$ with different values.
During the implementation, the stopping criterion in all algorithms is as $\| F_k\|\leq \epsilon=10^{-6}$, or the number of iterations exceed 2,000.

\subsection{Large-Scale Nonlinear Equations}
In this part, we test these algorithms on the standard constrained nonlinear equations with different dimensions. The following test Problems 1-4 are respectively selected as the same as Problems 1, 3, 5 and 7 in \cite{WSLZ23}. The convex constraints of these problems are $\mathcal{C}=\mathbb{R}_{+}^{n}$ and the mapping $F$ is defined as
\[
F(x)=\left(f_{1}(x),f_{2}(x),\cdots,f_{n}(x)\right)^\top.
\]

\noindent\textbf{Problem 1.} \[f_{i}(x)=e^{x_i}-1,~i=1,2,\cdots,n.\]

\noindent\textbf{Problem 2.} \[f_{i}(x)=\ln(x_i+1)-\frac{x_i}{n},~i=1,2,\cdots,n.\]

\noindent\textbf{Problem 3.} \[f_{1}(x)=e^{x_1}-1,~f_{i}(x)=e^{x_i}+x_i-1,~i=2,3,\cdots,n.\]

\noindent\textbf{Problem 4.} \[f_{i}(x)= 2x_i -\sin (x_i),~i=1,2,\cdots,n.\]

To assess the effectiveness of these algorithms objectively, we conduct tests for each problem using initial points randomly generated from the interval $(0, 1)$. The numerical results, obtained from running each test 10 times with each algorithm, are presented in Table \ref{T1}, where `` P(n)/$\overline{\rm Iter}$/$\overline{\rm NF}$/$\overline{\rm Tcpu}$/$\overline{\|F^{\ast}\|}$'' stand for test problems (problem dimensions), average number of iterations, average number of evaluations of $F$, average CPU time in seconds, average final value of $\|F_k\|$ when the program is stopped, respectively. Table \ref{T1} shows that the AA variants of three DFPM are all superior (in terms of $\overline{\rm Iter}$, $\overline{\rm NF}$ and $\overline{\|F^{\ast}\|}$) to their originals for these chosen set of test problems, which also confirms the encouraging capability of AA for DFPM. In contrast, $\overline{\rm Tcpu}$ deteriorates in certain tests as a result of AA having to solve an extra optimization problem in each iteration. 

\begin{figure}[ht]
	\subfigure[Performance profiles on $\overline{\rm Iter}$\label{F1a}]{
		\begin{minipage}{0.493\linewidth}
			\includegraphics[width=1\textwidth]{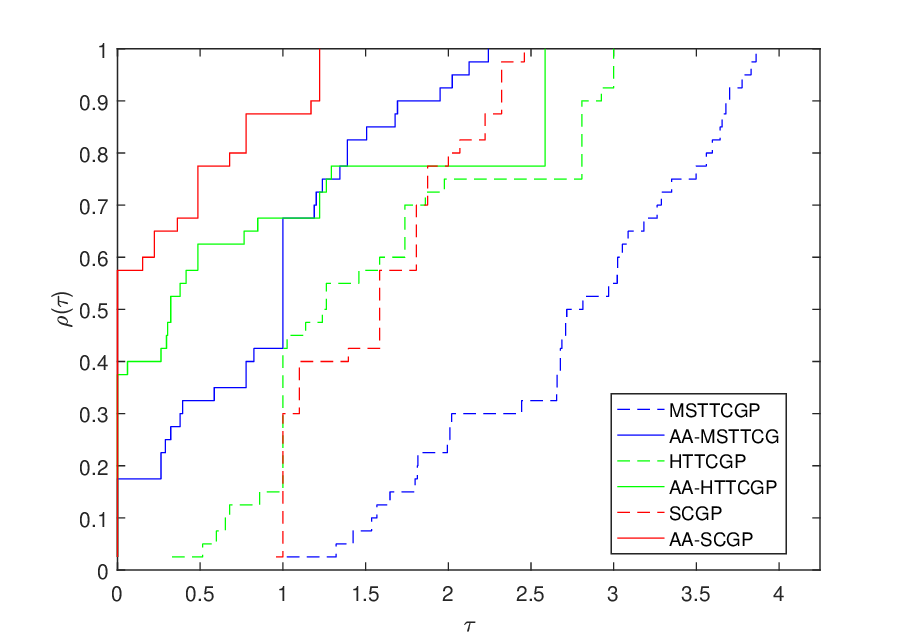}
	\end{minipage}}
	\subfigure[Performance profiles on $\overline{\rm NF}$\label{F1b}]{
		\begin{minipage}{0.493\linewidth}
			\includegraphics[width=1\textwidth]{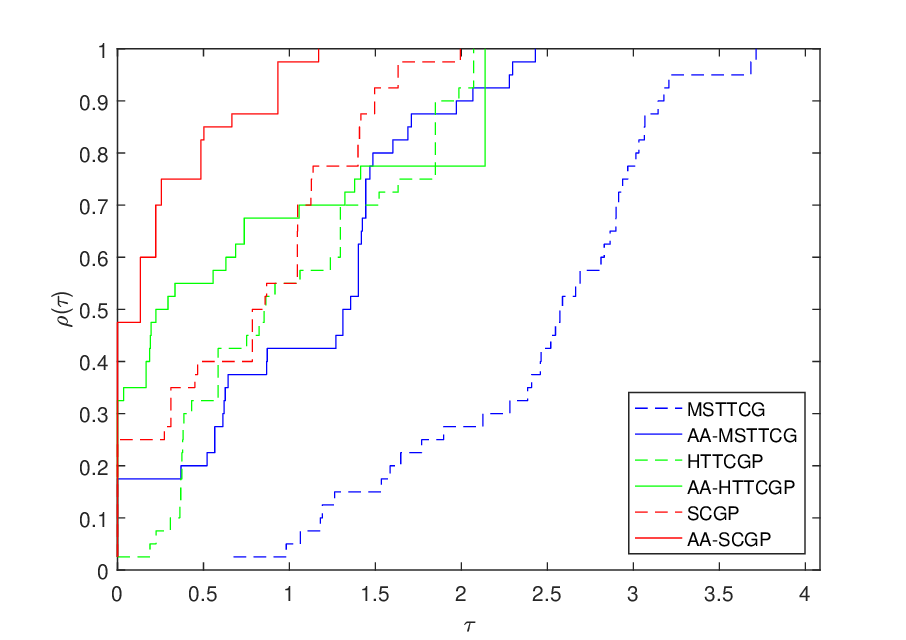}
	\end{minipage}}
	\caption{Performance profiles of these methods for constrained nonlinear equations}\label{F1}
\end{figure}

Moreover, we use the performance profiles \cite{DM02} to visually compare the performance of these methods, as illustrated in Figures \ref{F1a} and \ref{F1b}, which intuitively describe $\overline{\rm Iter}$ and $\overline{\rm NF}$, respectively.  
The performance profiles $\rho(\tau)$ show the probability that a solver is within a certain factor $\tau$ of the best possible performance. In short, the higher the curve, the better the method.
It is very clear from Figure \ref{F1} that the acceleration process is efficient in its purpose of accelerating DFPM.

\begin{landscape}
	\begin{table}
		\centering
		\caption{Numerical results on Problem 1-4 with random initial points}\label{T1}
		\resizebox{1.65\textwidth}{!}{
			\begin{tabular}{ccccccc}
				\hline
				\multirow{2}*{P(n)} & MSTTCGP & AA-MSTTCGP & HTTCGP & AA-HTTCGP & SCGP & AA-SCGP  \\
				& $\overline{\rm Iter}$/$\overline{\rm NF}$/$\overline{\rm Tcpu}$/$\overline{\|F^{\ast}\|}$
				& $\overline{\rm Iter}$/$\overline{\rm NF}$/$\overline{\rm Tcpu}$/$\overline{\|F^{\ast}\|}$
				& $\overline{\rm Iter}$/$\overline{\rm NF}$/$\overline{\rm Tcpu}$/$\overline{\|F^{\ast}\|}$
				& $\overline{\rm Iter}$/$\overline{\rm NF}$/$\overline{\rm Tcpu}$/$\overline{\|F^{\ast}\|}$
				& $\overline{\rm Iter}$/$\overline{\rm NF}$/$\overline{\rm Tcpu}$/$\overline{\|F^{\ast}\|}$
				& $\overline{\rm Iter}$/$\overline{\rm NF}$/$\overline{\rm Tcpu}$/$\overline{\|F^{\ast}\|}$\\
				\hline
				1(10000)&38.4/113.8/0.051/5.27e-07&5.0/23.0/0.048/{\bf0.00e+00}&11.0/29.1/{\bf0.015}/3.17e-07&{\bf4.0}/{\bf14.0}/0.029/{\bf0.00e+00}&14.0/29.0/0.020/6.56e-07&9.0/29.0/0.053/{\bf0.00e+00}\\
				1(30000)&44.8/135.0/0.135/5.78e-07&9.1/32.9/0.077/{\bf0.00e+00}&8.8/22.7/{\bf0.024}/3.37e-07&{\bf7.0}/{\bf22.0}/0.056/{\bf0.00e+00}&15.0/31.0/0.038/3.34e-07&9.0/29.0/0.076/{\bf0.00e+00}\\
				1(50000)&43.6/131.2/0.216/6.38e-07&9.7/30.7/0.104/{\bf0.00e+00}&11.8/31.0/0.051/4.13e-07&7.0/22.0/0.078/{\bf0.00e+00}&15.0/31.0/0.068/4.31e-07&{\bf3.0}/{\bf9.0}/{\bf0.033}/{\bf0.00e+00}\\
				1(80000)&42.7/128.4/0.298/4.62e-07&12.2/39.3/0.195/{\bf0.00e+00}&10.9/28.7/0.076/3.93e-07&7.0/21.0/0.114/1.99e-09&15.0/31.0/0.100/5.45e-07&{\bf3.0}/{\bf9.0}/{\bf0.040}/{\bf0.00e+00}\\
				1(100000)&40.6/125.4/0.358/5.48e-07&13.1/39.9/0.250/{\bf0.00e+00}&11.8/31.3/0.094/2.44e-07&7.0/21.0/0.139/8.52e-10&15.0/31.0/0.123/6.08e-07&{\bf5.0}/{\bf13.0}/{\bf0.092}/{\bf0.00e+00}\\
				1(120000)&45.4/138.5/0.501/4.22e-07&13.1/39.4/0.275/{\bf0.00e+00}&11.0/28.3/{\bf0.110}/2.60e-07&7.0/21.0/0.159/3.46e-10&15.0/31.0/0.144/6.67e-07&{\bf5.0}/{\bf13.0}/0.111/{\bf0.00e+00}\\
				1(150000)&47.3/142.4/0.750/6.23e-07&13.0/39.0/0.352/{\bf0.00e+00}&11.6/30.7/0.178/3.79e-07&7.0/21.0/0.198/4.15e-14&15.0/31.0/0.247/7.47e-07&{\bf5.7}/{\bf14.4}/{\bf0.165}/{\bf0.00e+00}\\
				1(180000)&44.2/130.6/0.802/5.02e-07&12.4/37.2/0.398/{\bf0.00e+00}&12.7/33.7/{\bf0.230}/3.19e-07&{\bf7.0}/{\bf21.0}/0.243/2.76e-11&15.0/31.0/0.280/8.17e-07&12.0/27.0/0.452/{\bf0.00e+00}\\
				1(200000)&46.0/137.9/0.924/4.91e-07&12.0/36.0/0.429/{\bf0.00e+00}&11.2/29.2/{\bf0.223}/3.85e-07&{\bf7.0}/{\bf21.0}/0.272/6.60e-15&15.0/31.0/0.327/8.61e-07&12.0/27.0/0.481/{\bf0.00e+00}\\
				1(250000)&44.1/132.7/1.124/5.09e-07&12.0/36.0/0.514/{\bf0.00e+00}&11.0/28.5/{\bf0.267}/2.31e-07&{\bf7.0}/{\bf21.0}/0.324/2.58e-14&15.0/31.0/0.396/9.63e-07&12.0/27.0/0.581/{\bf0.00e+00}\\
				2(10000)&25.3/81.8/0.024/3.10e-07&8.8/30.3/0.051/{\bf0.00e+00}&10.9/29.7/{\bf0.009}/5.77e-07&{\bf7.2}/{\bf19.4}/0.040/{\bf0.00e+00}&14.0/29.0/0.012/3.22e-07&8.0/19.0/0.045/{\bf0.00e+00}\\
				2(30000)&28.4/91.6/0.049/6.15e-07&9.2/31.4/0.064/{\bf0.00e+00}&10.0/27.0/{\bf0.016}/{\bf0.00e+00}&7.3/18.6/0.047/2.98e-07&14.0/29.0/0.022/5.61e-07&{\bf7.0}/{\bf17.0}/0.050/{\bf0.00e+00}\\
				2(50000)&24.4/82.0/0.067/3.42e-07&14.2/49.1/0.136/{\bf0.00e+00}&10.0/27.0/0.024/{\bf0.00e+00}&{\bf3.0}/{\bf10.0}/{\bf0.022}/{\bf0.00e+00}&14.0/29.0/0.032/7.23e-07&7.0/17.0/0.060/{\bf0.00e+00}\\
				2(80000)&21.1/71.0/0.091/3.96e-07&11.6/46.5/0.157/{\bf0.00e+00}&10.0/27.0/0.037/{\bf0.00e+00}&{\bf3.0}/{\bf10.0}/{\bf0.028}/{\bf0.00e+00}&14.0/29.0/0.050/9.14e-07&7.0/17.0/0.085/{\bf0.00e+00}\\
				2(100000)&29.3/99.5/0.162/2.19e-07&13.1/43.9/0.216/{\bf0.00e+00}&10.0/27.0/0.049/{\bf0.00e+00}&{\bf3.0}/{\bf10.0}/{\bf0.038}/{\bf0.00e+00}0&15.0/31.0/0.069/3.09e-07&7.0/17.0/0.108/{\bf0.00e+00}\\
				2(120000)&23.5/80.1/0.159/1.49e-07&9.6/37.4/0.183/{\bf0.00e+00}&10.0/27.0/0.057/{\bf0.00e+00}&{\bf3.0}/{\bf10.0}/{\bf0.043}/{\bf0.00e+00}&15.0/31.0/0.083/3.38e-07&7.0/17.0/0.119/{\bf0.00e+00}\\
				2(150000)&27.2/87.5/0.339/2.94e-07&12.7/46.4/0.366/{\bf0.00e+00}&10.0/27.0/{\bf0.118}/{\bf0.00e+00}&{\bf5.0}/{\bf16.0}/0.125/{\bf0.00e+00}&15.0/31.0/0.184/3.79e-07&7.0/17.0/0.192/{\bf0.00e+00}\\
				2(180000)&32.0/103.2/0.475/4.15e-07&11.5/40.9/0.380/{\bf0.00e+00}&10.0/27.0/{\bf0.140}/{\bf0.00e+00}&{\bf5.0}/{\bf16.0}/0.147/{\bf0.00e+00}&15.0/31.0/0.220/4.15e-07&7.0/17.0/0.222/{\bf0.00e+00}\\
				2(200000)&32.2/105.2/0.530/4.73e-07&11.8/40.4/0.428/{\bf0.00e+00}&10.0/27.0/{\bf0.153}/{\bf0.00e+00}&{\bf5.0}/{\bf16.0}/0.159/{\bf0.00e+00}&15.0/31.0/0.241/4.37e-07&7.0/17.0/0.244/{\bf0.00e+00}\\
				2(250000)&32.7/108.2/0.676/4.19e-07&14.2/48.3/0.631/2.44e-08&10.0/27.0/{\bf0.190}/{\bf0.00e+00}&{\bf5.0}/{\bf16.0}/0.194/{\bf0.00e+00}&15.0/31.0/0.293/4.89e-07&7.0/17.0/0.293/{\bf0.00e+00}\\
				3(10000)&17.4/69.0/0.030/6.67e-08&6.0/31.0/0.040/{\bf0.00e+00}&12.0/40.7/{\bf0.018}/1.78e-07&6.0/25.0/0.038/3.70e-14&20.0/61.0/0.030/7.03e-07&{\bf5.0}/{\bf19.0}/0.030/{\bf0.00e+00}\\
				3(30000)&17.6/69.5/0.066/7.75e-08&6.0/30.0/0.057/{\bf0.00e+00}&12.0/40.4/{\bf0.041}/1.56e-07&{\bf5.0}/{\bf21.0}/0.044/9.24e-10&21.0/64.0/0.075/4.90e-07&8.0/28.0/0.075/{\bf0.00e+00}\\
				3(50000)&17.8/70.8/0.109/2.09e-15&{\bf6.0}/{\bf30.0}/0.078/{\bf0.00e+00}&12.0/40.5/{\bf0.058}/1.63e-07&10.8/43.1/0.145/2.04e-08&21.0/64.0/0.111/6.30e-07&{\bf6.0}/{\bf30.0}/0.087/{\bf0.00e+00}\\
				3(80000)&18.8/74.4/0.177/1.63e-07&{\bf6.0}/{\bf30.0}/0.106/{\bf0.00e+00}&12.0/40.2/{\bf0.087}/2.42e-07&10.2/41.8/0.192/3.18e-08&21.0/64.0/0.167/7.99e-07&{\bf6.0}/{\bf30.0}/0.122/{\bf0.00e+00}\\
				3(100000)&12.2/49.5/0.141/8.46e-08&{\bf6.0}/{\bf30.0}/0.126/{\bf0.00e+00}&12.0/40.4/{\bf0.112}/2.69e-07&8.0/34.2/0.184/{\bf0.00e+00}&21.0/64.0/0.200/8.91e-07&{\bf6.0}/{\bf30.0}/0.152/{\bf0.00e+00}\\
				3(120000)&23.9/97.1/0.333/1.73e-07&{\bf6.0}/{\bf30.0}/0.153/{\bf0.00e+00}&12.0/40.0/{\bf0.137}/3.26e-07&7.5/32.2/0.193/{\bf0.00e+00}&21.0/64.0/0.255/9.76e-07&{\bf6.0}/{\bf30.0}/0.175/{\bf0.00e+00}\\
				3(150000)&16.1/64.8/0.344/8.86e-08&{\bf6.0}/{\bf30.0}/0.211/{\bf0.00e+00}&12.0/40.3/{\bf0.210}/3.02e-07&7.4/31.7/0.261/6.99e-09&22.0/67.6/0.421/5.59e-07&{\bf6.0}/{\bf30.0}/0.242/{\bf0.00e+00}\\
				3(180000)&15.0/61.2/0.359/4.86e-08&{\bf6.0}/{\bf30.0}/{\bf0.245}/{\bf0.00e+00}&12.0/40.2/0.254/3.68e-07&7.5/31.7/0.303/{\bf0.00e+00}&22.0/67.7/0.491/5.19e-07&7.0/33.0/0.334/{\bf0.00e+00}\\
				3(200000)&17.4/70.2/0.464/9.00e-08&{\bf6.0}/{\bf30.0}/{\bf0.268}/{\bf0.00e+00}&12.0/40.0/0.288/3.76e-07&7.8/32.8/0.325/1.72e-11&22.0/68.2/0.567/5.97e-07&7.0/33.0/0.350/{\bf0.00e+00}\\
				3(250000)&16.1/63.3/0.538/1.60e-08&{\bf6.0}/{\bf30.0}/0.355/{\bf0.00e+00}&12.0/40.1/{\bf0.349}/4.54e-07&9.8/40.6/0.531/{\bf0.00e+00}&22.0/67.7/0.694/6.13e-07&{\bf4.0}/{\bf16.0}/{\bf0.214}/{\bf0.00e+00}\\
				4(10000)&8.5/27.5/0.006/4.99e-07&2.0/11.4/0.006/{\bf0.00e+00}&7.0/18.0/0.004/2.62e-07&2.4/9.1/0.008/{\bf0.00e+00}&2.0/5.0/0.003/{\bf0.00e+00}&{\bf1.0}/{\bf5.0}/{\bf0.002}/{\bf0.00e+00}\\
				4(30000)&10.2/31.7/0.010/4.81e-07&2.0/12.6/0.008/{\bf0.00e+00}&7.0/18.0/0.007/4.50e-07&6.0/18.0/0.035/3.36e-09&2.0/5.0/0.005/{\bf0.00e+00}&{\bf1.0}/{\bf5.0}/{\bf0.002}/{\bf0.00e+00}\\
				4(50000)&11.3/35.1/0.016/4.08e-07&2.0/12.2/0.010/{\bf0.00e+00}&7.0/18.0/0.009/5.61e-07&6.0/18.0/0.042/2.95e-10&2.0/5.0/0.007/{\bf0.00e+00}&{\bf1.0}/{\bf5.0}/{\bf0.003}/{\bf0.00e+00}\\
				4(80000)&12.1/37.7/0.027/4.03e-07&2.0/11.4/0.015/{\bf0.00e+00}&7.0/18.0/0.015/7.25e-07&6.0/18.0/0.054/2.96e-10&2.0/5.0/0.011/{\bf0.00e+00}&{\bf1.0}/{\bf5.0}/{\bf0.005}/{\bf0.00e+00}\\
				4(100000)&12.6/39.1/0.035/5.64e-07&2.0/12.2/0.019/{\bf0.00e+00}&7.0/18.0/0.019/8.38e-07&6.0/18.0/0.068/2.87e-10&2.0/5.0/0.014/{\bf0.00e+00}&{\bf1.0}/{\bf5.0}/{\bf0.007}/{\bf0.00e+00}\\
				4(120000)&12.8/40.4/0.045/5.24e-07&2.0/12.6/0.024/{\bf0.00e+00}&7.0/18.0/0.023/9.00e-07&6.0/18.0/0.079/2.81e-10&2.0/5.0/0.017/{\bf0.00e+00}&{\bf1.0}/{\bf5.0}/{\bf0.009}/{\bf0.00e+00}\\
				4(150000)&11.8/37.3/0.089/4.04e-07&2.0/11.8/0.037/{\bf0.00e+00}&7.6/19.8/0.055/5.14e-07&6.0/18.0/0.116/5.67e-10&2.0/5.0/0.030/{\bf0.00e+00}&{\bf1.0}/{\bf5.0}/{\bf0.017}/{\bf0.00e+00}\\
				4(180000)&12.5/38.3/0.110/3.37e-07&2.0/12.2/0.046/{\bf0.00e+00}&8.0/21.0/0.070/2.25e-07&6.0/18.0/0.140/7.78e-10&2.0/5.0/0.036/{\bf0.00e+00}&{\bf1.0}/{\bf5.0}/{\bf0.020}/{\bf0.00e+00}\\
				4(200000)&13.0/40.9/0.128/2.93e-07&2.0/12.6/0.049/{\bf0.00e+00}&8.0/21.0/0.074/2.34e-07&6.0/18.0/0.146/9.32e-10&2.0/5.0/0.040/{\bf0.00e+00}&{\bf1.0}/{\bf5.0}/{\bf0.023}/{\bf0.00e+00}\\
				4(250000)&13.7/44.2/0.173/3.57e-07&2.0/12.2/0.060/{\bf0.00e+00}&8.0/21.0/0.093/2.64e-07&6.0/18.0/0.178/1.03e-09&2.0/5.0/0.050/{\bf0.00e+00}&{\bf1.0}/{\bf5.0}/{\bf0.027}/{\bf0.00e+00}\\
				\hline
			\end{tabular}
		}
\end{table}\end{landscape}

\begin{figure}[ht]
		\subfigure[Performance profiles on $\overline{\rm Iter}$\label{F3a}]{
			\begin{minipage}{0.493\linewidth}
				\includegraphics[width=1\textwidth]{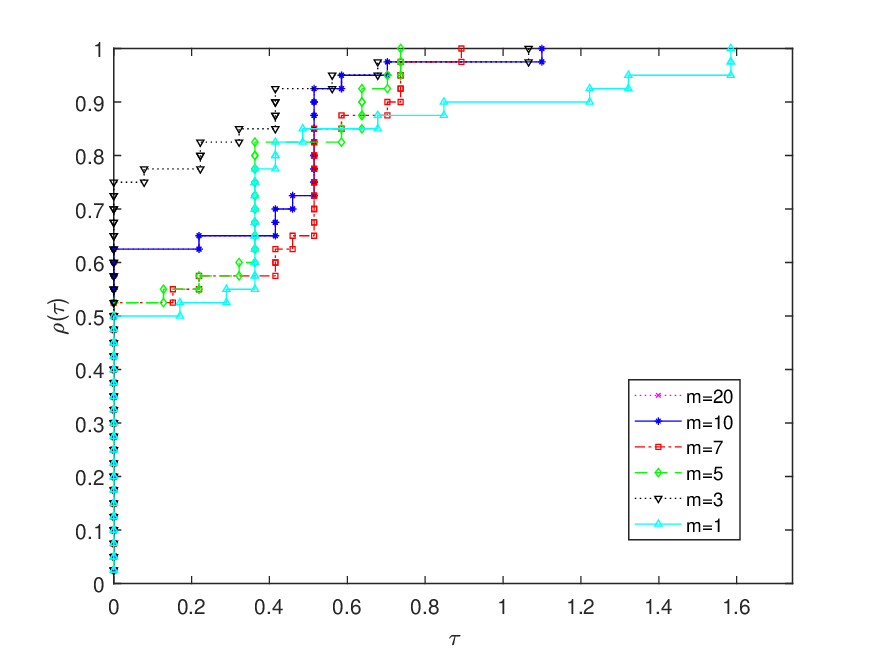}
		\end{minipage}}
		\subfigure[Performance profiles on $\overline{\rm NF}$\label{F3b}]{
			\begin{minipage}{0.493\linewidth}
				\includegraphics[width=1\textwidth]{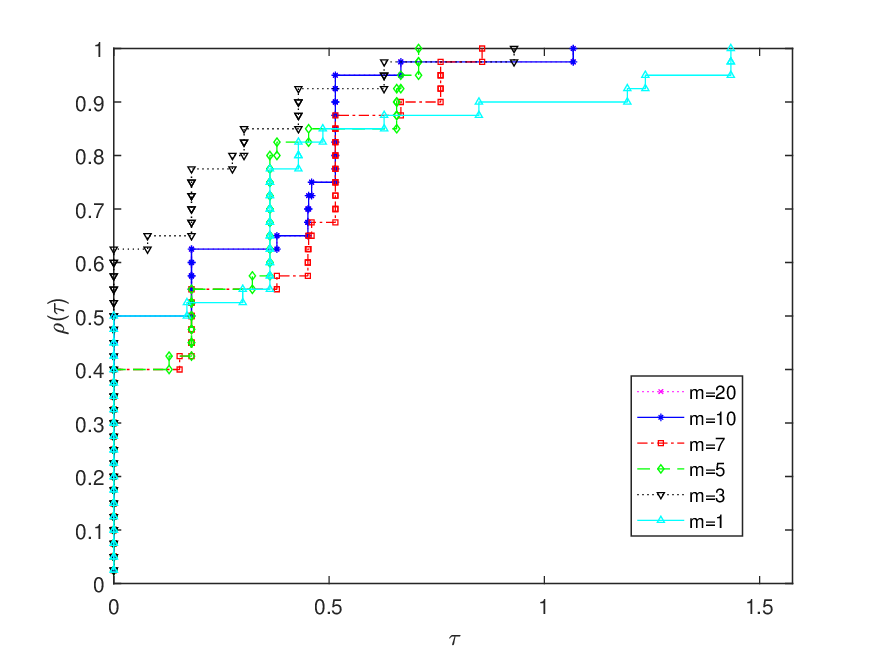}
		\end{minipage}}
		\caption{The effect of $m$ for AA-SCGP on solving constrained nonlinear equations}\label{F2}
\end{figure}

Next, we consider the impact of the choice of $m$. Leave other parameters unchanged. Performance profiles on $\overline{\rm Iter}$ and $\overline{\rm NF}$ for AA-SCGP with $m=1,3,5,7,10,20$ are plotted in Figure \ref{F2}, from which we can see that $m=3$ is the best. We favor the modest values of $m$. In nonlinear problems, the inclusion of unrepresentative older iterants may be detrimental, and large $m$ can cause numerical difficulties in acceleration.

\subsection{Regularized Decentralized Logistic Regression}
We consider a real-world application, regularized decentralized logistic regression, which is a classic example that is widely used \cite{JYTH22,Luo21,YJM23}.
\begin{equation}\label{RDLR}
	\min _{x \in \mathbb{R}^n} f(x):=\frac{1}{T} \sum_{i=1}^T \ln (1+\exp (-b_i a_i^{\top} x))+\frac{\tau}{2}\|x\|^2,
\end{equation}
where $\tau>0$ is a regularization parameter, $\frac{1}{T} \sum_{i=1}^T \ln (1+\exp (-b_i a_i^{\top} x))$ represents the logistic loss function, and the data pairs $(a_i, b_i) \in \mathbb{R}^n \times\{-1,1\}(i=1, \ldots, T)$ are taken from a given data set or distribution. It is easy to know that the objective function $f$ is strongly convex and has Lipschitz continuous gradient \cite{Luo21}. Hence $x^* \in \mathbb{R}^n$ is a unique optimal solution to \eqref{RDLR} if and only if it is a root of the following nonlinear equations \cite{JYTH22},
\begin{equation}\label{RDLRe}
	F(x):=\nabla f(x)=\frac{1}{T} \sum_{i=1}^T \frac{-b_i \exp (-b_i a_i^{\top} x) a_i}{1+\exp (-b_i a_i^{\top} x)}+\tau x=0.
\end{equation}
The problem \eqref{RDLRe} is strongly monotone and Lipschitz continuous, thus it satisfies the local error bound Assumption \ref{ass2}, and we can apply AA-DFPM to solve the above problem.
Considering that problem \eqref{RDLRe} is unconstrained, we can dispense with the nonnegative constraints in \eqref{comak}. Let $a^k_k=1-\sum\limits_{j=k-{m_k}}^{k-1} a^k_j$, so the least squares problem can be reformulated as
\begin{equation}\label{ULSP}
	\min\limits_{(a^k_{k-{m_k}},\dots,a^k_{k-1})^\top}\left\|r_k+\sum\limits_{j=k-{m_k}}^{k-1} a^k_j (r_j-r_k)\right\|^2.
\end{equation}
Our implementation solves above problem using QR decomposition. The QR decomposition of problem \eqref{ULSP} at iteration $k$ can be efficiently obtained from that of at iteration $k-1$ in $O(m_k n)$ \cite{GVL96}.

We exclusively focus on AA-SCGP, the top-performing method in our initial experiments, and compare it with two DFPM incorporating inertial acceleration: MITTCGP \cite{MJJYH23} and IHCGPM3 \cite{JYTH22}. The involved parameters for both MITTCGP and IHCGPM3 are set to their defaults, while the parameters used in AA-SCGP are taken from the experiment in last subsection. The test instances are sourced from the LIBSVM datasets\footnote{Datasets available at \href{https://www.csie.ntu.edu.tw/~cjlin/libsvmtools/datasets/}{https://www.csie.ntu.edu.tw/$\sim$cjlin/libsvmtools/datasets/}.} \cite{CL11} and the termination criterion of all three algorithms is the same as the first experiment.

\begin{table}[h]		
		\caption{The effect of $m$ for AA-SCGP on solving problem \ref{RDLRe}}\label{T2}
		\centering
		\resizebox{\textwidth}{!}{
			\begin{tabular}{cccc}
				\hline
				\multirow{2}*{Data sets}  & $m=1$ & $m=3$ & $m=5$ \\
				& Iter/NF/Tcpu/$\|F^{\ast}\|$
				& Iter/NF/Tcpu/$\|F^{\ast}\|$
				& Iter/NF/Tcpu/$\|F^{\ast}\|$\\
				\hline
				fourclass\_scale & 156.0/468.0/0.052/9.50e-07 & {\bf11.0/33.0/0.022}/2.78e-07 & 15.0/46.0/0.026/{\bf6.38e-08} \\ 
				liver-disorders & 313.0/980.0/{\bf0.029}/9.69e-07 & 255.0/790.0/0.048/9.54e-07 & {\bf141.0/452.0}/0.039/{\bf3.72e-07}\\ 
				phishing & ---& 115.0/343.0/3.790/9.47e-07 & {\bf55.0/163.0/1.840}/9.51e-07 \\ 
				w4a & 218.0/651.0/{\bf1.056}/9.88e-07 & 312.0/934.0/1.887/8.96e-07 & 903.0/2705.0/5.609/{\bf4.51e-07} \\ 
				w5a & {\bf211.0/630.0/1.350}/9.79e-07 & 980.0/2934.0/7.898/7.36e-07 & 1552.0/4649.0/12.785/9.39e-07\\ 
				w6a & {\bf207.0/618.0/2.995}/9.12e-07 & 600.0/1795.0/10.337/{\bf6.47e-07} & --- \\ 
				\hline
				\multirow{2}*{Data sets}  & $m=7$ & $m=10$ & $m=20$\\
				& Iter/NF/Tcpu/$\|F^{\ast}\|$
				& Iter/NF/Tcpu/$\|F^{\ast}\|$
				& Iter/NF/Tcpu/$\|F^{\ast}\|$\\
				\hline
				fourclass\_scale & 14.0/42.0/0.022/5.66e-07 & 17.0/51.0/0.023/6.77e-07 & 27.0/81.0/0.024/9.38e-07\\ 
				liver-disorders & 718.0/2196.0/0.110/6.10e-07 & 195.0/614.0/0.043/6.14e-07 & ---\\ 
				phishing & 152.0/451.0/4.579/4.69e-07 & 154.0/455.0/4.424/8.27e-07 & 424.0/1264.0/12.295/{\bf2.65e-07}\\ 
				w4a & ---& {\bf185.0/546.0}/1.223/9.65e-07 & 284.0/841.0/2.201/9.70e-07\\ 
				w5a & 988.0/2957.0/8.298/{\bf5.69e-07} & 1302.0/3897.0/11.513/8.29e-07 & 1893.0/5668.0/18.758/7.66e-07\\ 
				w6a & --- & 423.0/1264.0/7.540/9.00e-07 & 270.0/799.0/5.068/9.34e-07\\ 
				\hline
			\end{tabular}
	}
\end{table}

\begin{figure}[ht]
	\subfigure[fourclass\_scale: $m=3$\label{F2a}]{
		\begin{minipage}{0.493\linewidth}
			\includegraphics[width=0.95\textwidth]{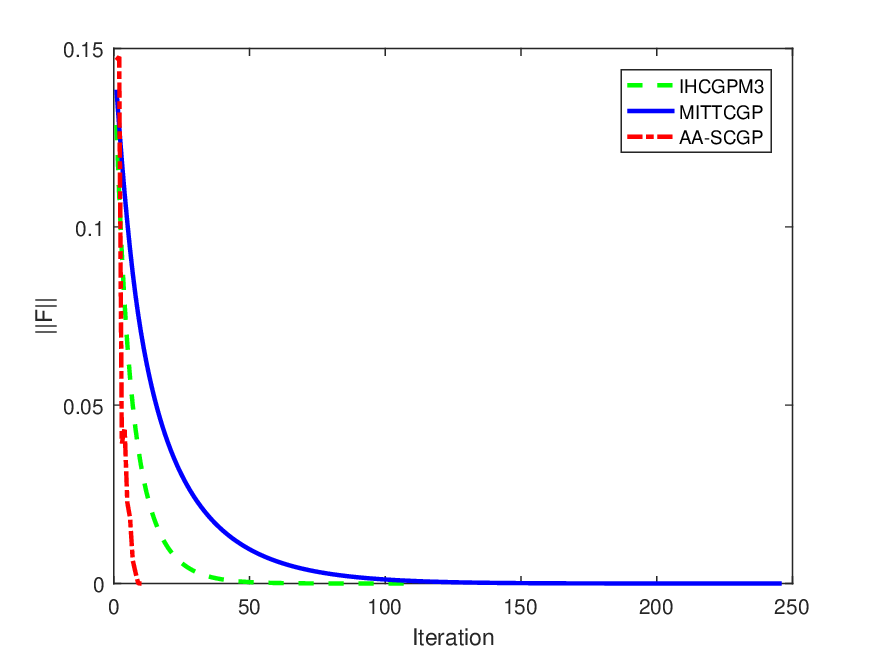}
	\end{minipage}}
	\subfigure[phishing: $m=3$\label{F2b}]{
		\begin{minipage}{0.493\linewidth}
			\includegraphics[width=0.95\textwidth]{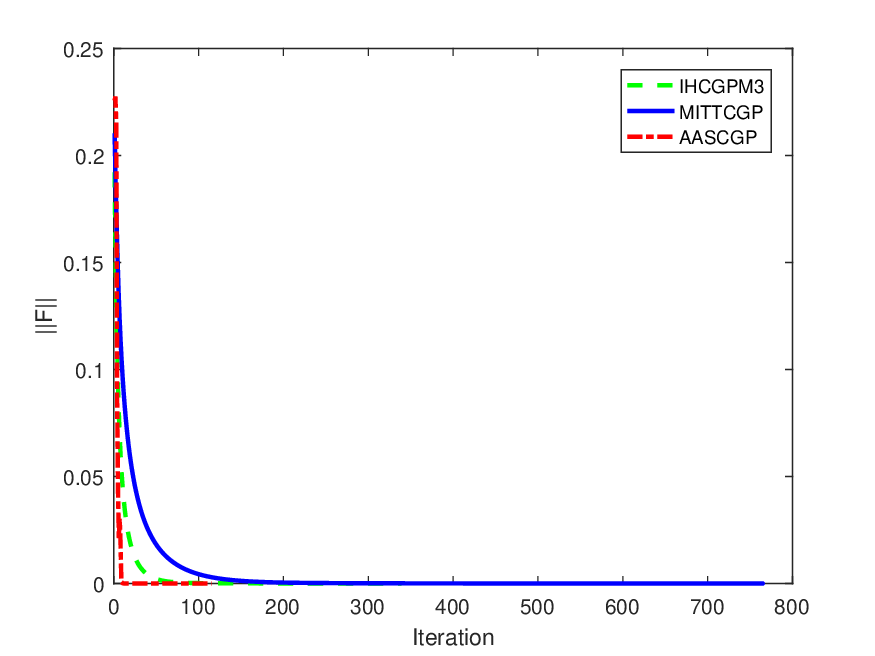}
	\end{minipage}}
	\subfigure[w5a: $m=3$\label{F2c}]{
		\begin{minipage}{0.493\linewidth}
			\includegraphics[width=0.95\textwidth]{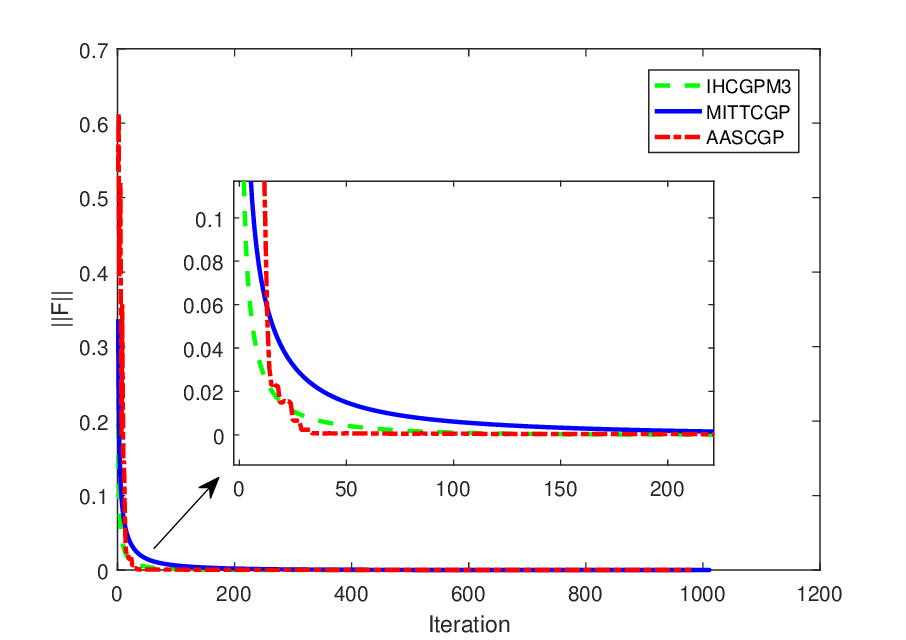}
	\end{minipage}}
	\subfigure[w6a: $m=3$\label{F2d}]{
		\begin{minipage}{0.493\linewidth}
			\includegraphics[width=0.95\textwidth]{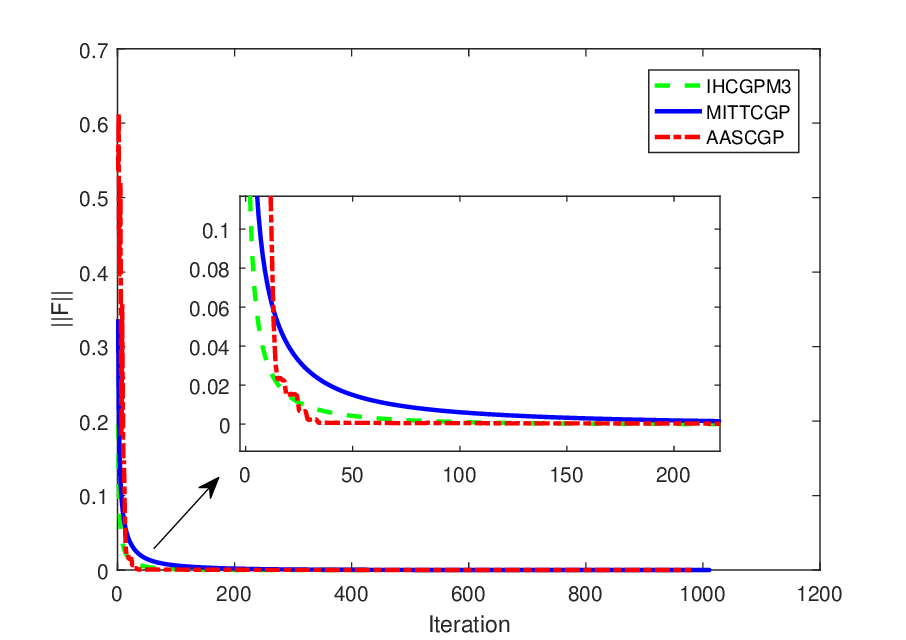}
	\end{minipage}}
		\subfigure[madelon.t $(T=600,n=500)$: $m=4$\label{F2e}]{
			\begin{minipage}{0.493\linewidth}
				\includegraphics[width=0.95\textwidth]{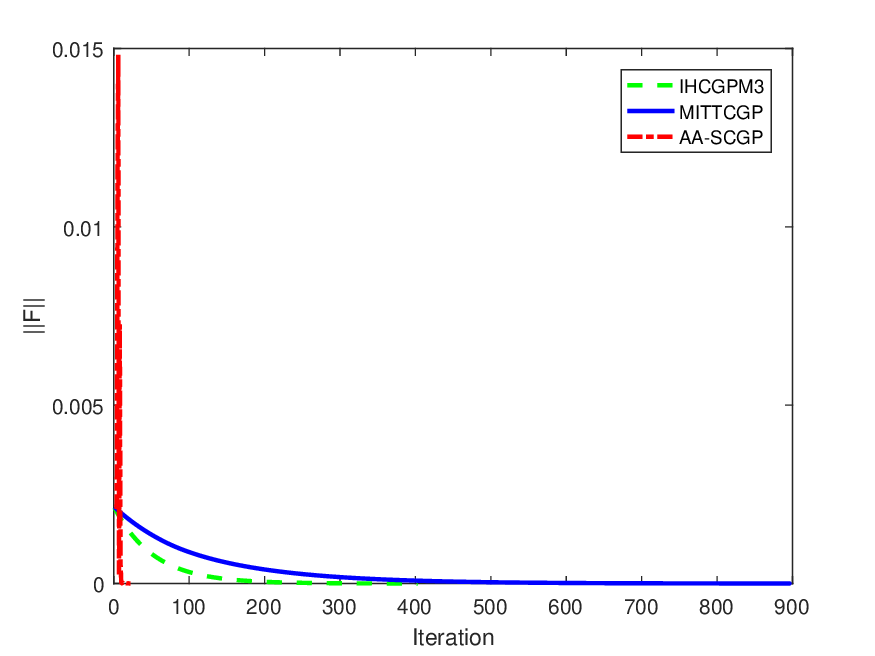}
		\end{minipage}}
		\subfigure[colon-cancer $(T=62,n=2000)$: $m=2$\label{F2f}]{
			\begin{minipage}{0.493\linewidth}
				\includegraphics[width=0.95\textwidth]{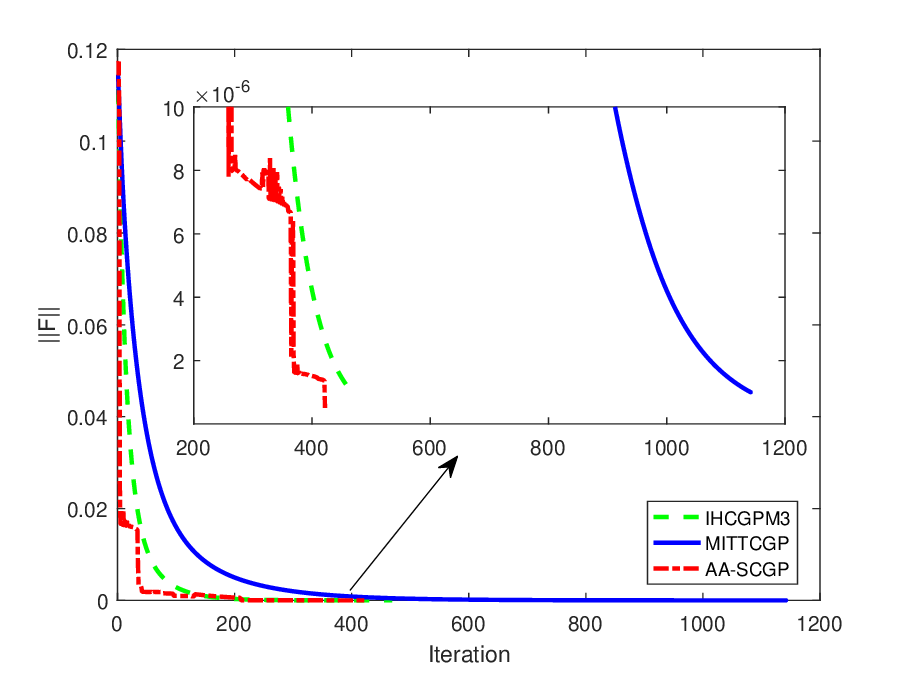}
		\end{minipage}}
		\caption{Change of $\|F(x)\|$ for problem \eqref{RDLRe} with initial point $(0,\dots,0)^\top$ and $\tau=0.01$}\label{F3}
\end{figure}

	First, we test the performance of AA-SCGP with different $m$. Set the origin as the initial point and $\tau = 0.01$. The results for AA-SCGP with $m = 1,3,5,7,10,20$ is showed in Table \ref{T2}, where Iter/NF/Tcpu/$\|F^{\ast}\|$ stand for number of iterations, number of evaluations of $F$, CPU time in seconds, final value of $\|F_k\|$ when the program is stopped, whereas ``---'' indicates a failure. 
	From Table \ref{T2}, we can see that the number of Iter/NF/Tcpu/$\|F^{\ast}\|$ does not decrease monotonically as $m$ increases. The point of diminishing returns is problem dependent and is perhaps best chosen by preliminary experiments for given problems. Whether the dynamic selection approaches can improve the convergence of AA-DFPM is an interesting topic for further research.

Next, we compare the performance of AA-SCGP with MITTCGP and IHCGPM3. Figure \ref{F3} displays the results of six test instances solved by three methods with a fixed initial point $(0,\dots,0)^\top$ and $\tau=0.01$. Here madelon.t and colon-cancer have been preprocessed and normalized. One can observe that our method outperforms the others. In particular, AA-SCGP performs better on different datasets fourclass\_scale, madelon.t and colon-cancer. The reason for this could be that these datasets have been scaled to [-1,1] or [0,1].

\begin{table}[h]		
		\caption{Numerical results for problem \eqref{RDLRe} with synthetic data}\label{T3}
		\centering
		\resizebox{\textwidth}{!}{
			\begin{tabular}{cccc}
				\hline
				\multirow{2}*{$(T,n)$} & IHCGPM3 &  MITTCGP & AA-SCGP($m=3$) \\
				& $\overline{\rm Iter}$/$\overline{\rm NF}$/$\overline{\rm Tcpu}$/$\overline{\|F^{\ast}\|}$
				& $\overline{\rm Iter}$/$\overline{\rm NF}$/$\overline{\rm Tcpu}$/$\overline{\|F^{\ast}\|}$
				& $\overline{\rm Iter}$($\overline{\rm NAA}$)/$\overline{\rm NF}$/$\overline{\rm Tcpu}$/$\overline{\|F^{\ast}\|}$\\
				\hline
				(500,1000) & 69.0/208.0/0.868/9.26e-07 & 171.0/514.0/2.070/9.56e-07 & {\bf25.4}(22.4)/{\bf74.4/0.618/6.28e-07}\\ 
				(1000,2000) & 71.0/214.0/3.311/8.68e-07 & 176.2/529.2/8.213/9.83e-07 & {\bf20.0}(17.0)/{\bf58.4/1.892/5.97e-07}\\ 
				(1500,3000) & 72.0/217.0/7.953/8.74e-07 & 179.8/538.8/19.857/9.92e-07 & {\bf20.6}(17.6)/{\bf60.6/4.522/6.48e-07}\\ 
				(2000,4000) & 72.6/217.6/13.856/9.78e-07 & 182.0/545.8/34.848/9.65e-07 & {\bf21.8}(18.8)/{\bf64.0/8.272/7.93e-07}\\ 
				(2500,5000) & 73.0/220.0/20.348/9.12e-07 & 183.0/550.0/50.693/9.82e-07 & {\bf17.2}(14.2)/{\bf50.6/9.725/7.36e-07}\\
				(5000,10000) & 75.0/224.8/73.730/9.42e-07 & 188.2/565.2/186.117/9.90e-07 & {\bf9.4}(6.4)/{\bf26.4/17.793/7.93e-07}\\ 
				(7500,15000) & 76.0/227.4/168.059/9.67e-07 & 191.0/574.0/423.680/9.92e-07 & {\bf13.2}(10.2)/{\bf37.6/59.334/7.89e-07}\\ 
				(10000,20000) & 76.0/229.0/276.792/9.79e-07 & 193.4/580.4/701.337/9.91e-07 & {\bf6.0}(2.0)/{\bf16.0/43.946/6.01e-07}\\ 
				(12500,25000) & 77.0/232.0/433.127/8.87e-07 & 195.0/586.0/1088.880/9.62e-07 & {\bf6.0}(2.0)/{\bf16.0/140.539/5.34e-07}\\
				\hline
			\end{tabular}
	}
\end{table}

Set $\tau = 0.01$ for real data and $\tau = 0.1$ for synthetic data. We use the MATLAB script ``{\rm 2*(rand(n,1)-0.5)}'' to generate the random initial point and run the same test 5 times for each test instance. Tables \ref{T3} and \ref{T4} show the numerical results, where the additional item $\overline{\rm NAA}$ indicates the average number of AA, and the other items are the same as in Table \ref{T1}. From Tables \ref{T3} and \ref{T4}, we can see that for most test instances, AA-SCGP outperforms two inertial methods in terms of $\overline{\rm NF}$ and $\overline{\rm Tcpu}$. In addition, the quality of the solutions obtained by AA-SCGP is better than that of the others. This benefits from the fact that AA-SCGP accelerates frequently during its iteration, in which the proportion of AA is 79.2\% for synthetic data and 55.6\% for real data. The numerical results also show that the improvement of AA-SCGP increases with $n$. These facts further illustrate that DFPM integrated with Anderson acceleration is valid and promising.

\begin{table}[h]		
		\caption{Numerical results for problem \eqref{RDLRe} with real data}\label{T4}
		\centering
		\resizebox{\textwidth}{!}{
			\begin{tabular}{ccccc}
				\hline
				\multirow{2}*{Data sets} & \multirow{2}*{$(T,n)$} & IHCGPM3 &  MITTCGP & AA-SCGP($m=3$) \\
				&
				& $\overline{\rm Iter}$/$\overline{\rm NF}$/$\overline{\rm Tcpu}$/$\overline{\|F^{\ast}\|}$
				& $\overline{\rm Iter}$/$\overline{\rm NF}$/$\overline{\rm Tcpu}$/$\overline{\|F^{\ast}\|}$
				& $\overline{\rm Iter}$($\overline{\rm NAA}$)/$\overline{\rm NF}$/$\overline{\rm Tcpu}$/$\overline{\|F^{\ast}\|}$\\
				\hline
				fourclass\_scale & (862,2) & 107.6/323.8/0.017/9.60e-07 & 246.6/740.4/0.038/9.82e-07 & {\bf11.4}(10.0)/{\bf33.8}/{\bf0.011}/{\bf2.18e-07}\\
				liver-disorders & (145,5) & {\bf492.2}/{\bf1520.6}/{\bf0.023}/9.93e-07 & 1362.6/4136.8/0.065/9.98e-07 & 760.0(462.4)/2548.8/0.088/{\bf7.70e-07}\\
				phishing & (11055,68) & {\bf561.8}/1685.6/14.016/9.93e-07 & 1394.0/4182.6/34.768/9.97e-07 & 565.0(243.4)/{\bf1374.4}/{\bf13.651}/{\bf5.14e-07}\\
				w4a & (7366,300) & {\bf584.6}/1754.8/{\bf2.815}/9.92e-07 & 1461.4/4384.4/7.050/9.97e-07 & 628.0(309.0)/{\bf1566.2}/3.104/{\bf9.11e-07}\\
				w5a & (9888,300) & 584.0/1752.6/3.714/9.95e-07 & 1459.8/4380.0/9.226/9.98e-07 & {\bf520.0}(87.0)/{\bf1128.0}/{\bf3.039}/{\bf9.53e-07}\\
				w6a & (17188,300) & {\bf584.4}/{\bf 1753.0}/{\bf8.294}/9.95e-07 & 1460.2/4381.6/20.713/9.96e-07 & 757.2(575.8)/2091.2/12.007/{\bf7.83e-07}\\
				\hline
			\end{tabular}
	}
\end{table}

\section{Conclusions}\label{Con}
In this paper, we developed a novel algorithm of using Anderson acceleration (AA) for derivative-free projection method (DFPM) in solving convex-constrained monotone nonlinear equations. First, we reviewed the convergence of a general framework for DFPM, and then explored how AA can still be exploited with DFPM though it may not a fixed-point iteration. As a result, an acceleration algorithm (AA-DFPM) with slight modifications is proposed, and the global convergence of AA-DFPM is obtained with no additional assumptions. The convergence rate is further established based on some suitable conditions. The results on both preliminary numerical experiments and applications demonstrate the superior performance.
As a future research, we plan to investigate a novel DFPM for general nonlinear equations. Considering that the least squares problem is hard to solve for a large window size $m$, we also intend to explore a simplified AA whose coefficients are easy to calculate. Moreover, designing different acceleration weights for $x_j$ and $v_j$ is one interesting topic.

\begin{acknowledgements}
This work was supported by the National Key Research and Development Program (2020YFA0713504) and the National Natural Science Foundation of China (12471401).

The authors would like to thank the two anonymous referees for their detailed reviews and insightful suggestions.
\end{acknowledgements}

%References
%\bibliographystyle{spmpsci}
%\bibliography{bibliography}

\begin{thebibliography}{10}
	\providecommand{\url}[1]{{#1}}
	\providecommand{\urlprefix}{URL }
	\expandafter\ifx\csname urlstyle\endcsname\relax
	\providecommand{\doi}[1]{DOI~\discretionary{}{}{}#1}\else
	\providecommand{\doi}{DOI~\discretionary{}{}{}\begingroup
		\urlstyle{rm}\Url}\fi
	
	\bibitem{Alvarez00}
	Alvarez, F.: On the minimizing property of a second order dissipative system in
	{H}ilbert spaces.
	\newblock SIAM J. Control Optim. \textbf{38}(4), 1102--1119 (2000)
	
	\bibitem{AF23}
	Amini, K., Faramarzi, P.: Global convergence of a modified spectral three-term
	{CG} algorithm for nonconvex unconstrained optimization problems.
	\newblock J. Comput. Appl. Math. \textbf{417}, 114630 (2023)
	
	\bibitem{AK15}
	Amini, K., Kamandi, A.: A new line search strategy for finding separating
	hyperplane in projection-based methods.
	\newblock Numer. Algorithms \textbf{70}, 559--570 (2015)
	
	\bibitem{Anderson65}
	Anderson, D.G.: Iterative procedures for nonlinear integral equations.
	\newblock J. Assoc. Comput. Mach. \textbf{12}(4), 547--560 (1965)
	
	\bibitem{Anderson19}
	Anderson, D.G.: Comments on ``{A}nderson acceleration, mixing and
	extrapolation''.
	\newblock Numer. Algorithms \textbf{80}, 135--234 (2019)
	
	\bibitem{BC22}
	Bian, W., Chen, X.J.: Anderson acceleration for nonsmooth fixed point problems.
	\newblock SIAM J. Numer. Anal. \textbf{60}(5), 2565--2591 (2022)
	
	\bibitem{CGH14}
	Cai, X.J., Gu, G.Y., He, B.S.: On the ${O}(1/t)$ convergence rate of the
	projection and contraction methods for variational inequalities with
	{L}ipschitz continuous monotone operators.
	\newblock Comput. Optim. Appl. \textbf{57}, 339--363 (2014)
	
	\bibitem{CL11}
	Chang, C.C., Lin, C.J.: {LIBSVM}: a library for support vector machines.
	\newblock ACM Trans. Intell. Syst. Technol. \textbf{2}(3), 1--27 (2011)
	
	\bibitem{CMY14}
	Chen, C.H., Ma, S.Q., Yang, J.F.: A general inertial proximal point algorithm
	for mixed variational inequality problem.
	\newblock SIAM J. Optim. \textbf{25}(4), 2120--2142 (2015)
	
	\bibitem{CZ14}
	Chorowski, J., Zurada, J.M.: Learning understandable neural networks with
	nonnegative weight constraints.
	\newblock IEEE Trans. Neural Netw. Learn. Syst. \textbf{26}(1), 62--69 (2014)
	
	\bibitem{DM02}
	Dolan, E.D., Mor\'{e}, J.J.: Benchmarking optimization software with
	performance profiles.
	\newblock Math. Program. \textbf{91}, 201--213 (2002)
	
	\bibitem{FS09}
	Fang, H.R., Saad, Y.: Two classes of multisecant methods for nonlinear
	acceleration.
	\newblock Numer. Linear Algebra Appl. \textbf{16}(3), 197--221 (2009)
	
	\bibitem{GLZ23}
	Garner, C., Lerman, G., Zhang, T.: Improved convergence rates of {A}nderson
	acceleration for a large class of fixed-point iterations (2023).
	\newblock Preprint at \url{https://arxiv.org/abs/2311.02490}
	
	\bibitem{GVL96}
	Golub, G.H., Van~Loan, C.: Matrix {C}omputations.
	\newblock The {J}ohns {H}opkins {U}niversity {P}ress, Baltimore (2013)
	
	\bibitem{GM23}
	Goncalves, M.L.N., Menezes, T.C.: A framework for convex-constrained monotone
	nonlinear equations and its special cases.
	\newblock Comp. Appl. Math. \textbf{42}, 306 (2023)
	
	\bibitem{HV19}
	Henderson, N.C., Varadhan, R.: Damped {A}nderson acceleration with restarts and
	monotonicity control for accelerating {EM} and {EM}-like algorithms.
	\newblock J. Comput. Graph. Stat. \textbf{28}(4), 834--846 (2019)
	
	\bibitem{IKRPA22}
	Ibrahim, A.H., Kumam, P., Rapaji\'{c}, S., Papp, Z., Abubakar, A.B.:
	Approximation methods with inertial term for large-scale nonlinear monotone
	equations.
	\newblock Appl. Numer. Math. \textbf{181}, 417--435 (2022)
	
	\bibitem{JYTH22}
	Jian, J.B., Yin, J.H., Tang, C.M., Han, D.L.: A family of inertial
	derivative-free projection methods for constrained nonlinear pseudo-monotone
	equations with applications.
	\newblock Comp. Appl. Math. \textbf{41}, 309 (2022)
	
	\bibitem{KSC02}
	Kudin, K.N., Scuseria, G.E., Cances, E.: A black-box self-consistent field convergence algorithm: One step closer.
	\newblock J. Chem. Phys. \textbf{116}(19), 8255--8261 (2002)
	
	\bibitem{LL11}
	Li, Q.N., Li, D.H.: A class of derivative-free methods for large-scale
	nonlinear monotone equations.
	\newblock IMA J. Numer. Anal. \textbf{31}(4), 1625--1635 (2011)
	
	\bibitem{LF19}
	Liu, J.K., Feng, Y.M.: A derivative-free iterative method for nonlinear
	monotone equations with convex constraints.
	\newblock Numer. Algorithms \textbf{82}, 245--262 (2019)
	
	\bibitem{LLXWT22}
	Liu, J.K., Lu, Z.L., Xu, J.L., Wu, S., Tu, Z.W.: An efficient projection-based
	algorithm without {L}ipschitz continuity for large-scale nonlinear
	pseudo-monotone equations.
	\newblock J. Comput. Appl. Math. \textbf{403}, 113822 (2022)
	
	\bibitem{Luo21}
	Luo, H.: Accelerated primal-dual methods for linearly constrained convex
	optimization problems (2021).
	\newblock Preprint at \url{https://arxiv.org/abs/2109.12604}
	
	\bibitem{MJJYH23}
	Ma, G.D., Jin, J.C., Jian, J.B., Yin, J.H., Han, D.L.: A modified inertial
	three-term conjugate gradient projection method for constrained nonlinear
	equations with applications in compressed sensing.
	\newblock Numer. Algorithms \textbf{92}, 1621--1653 (2023)
	
	\bibitem{MJ20}
	Mai, V., Johansson, M.: Anderson acceleration of proximal gradient methods.
	\newblock In: H.~Daum\'{e}~III, A.~Singh (eds.) Proceedings of the 37th
	International Conference on Machine Learning, \emph{Proceedings of Machine
		Learning Research}, vol. 119, pp. 6620--6629. PMLR (2020)
	
	\bibitem{MM87}
	Meintjes, K., Morgan, A.P.: A methodology for solving chemical equilibrium
	systems.
	\newblock Appl. Math. Comput. \textbf{22}(4), 333--361 (1987)
	
	\bibitem{Nesterov83}
	Nesterov, Y.: A method of solving a convex programming problem with convergence
	rate ${O}(1/k^2)$.
	\newblock Soviet Math. Dokl. \textbf{27}, 372--376 (1983)
	
	\bibitem{OL18}
	Ou, Y.G., Li, J.Y.: A new derivative-free {SCG}-type projection method for
	nonlinear monotone equations with convex constraints.
	\newblock J. Appl. Math. Comput. \textbf{56}, 195--216 (2018)
	
	\bibitem{OL23}
	Ou, Y.G., Li, L.: A unified convergence analysis of the derivative-free
	projection-based method for constrained nonlinear monotone equations.
	\newblock Numer. Algorithms \textbf{93}, 1639--1660 (2023)
	
	\bibitem{OPYZD20}
	Ouyang, W.Q., Peng, Y., Yao, Y.X., Zhang, J.Y., Deng, B.L.: Anderson
	acceleration for nonconvex {ADMM} based on {D}ouglas-{R}achford splitting.
	\newblock Comput. Graph. Forum \textbf{39}(5), 221--239 (2020)
	
	\bibitem{OTMD23}
	Ouyang, W.Q., Tao, J., Milzarek, A., Deng, B.L.: Nonmonotone globalization for
	{A}nderson acceleration via adaptive regularization.
	\newblock J. Sci. Comput. \textbf{96}, 5 (2023)
	
	\bibitem{SRX19}
	Pollock, S., Rebholz, L.G., Xiao, M.Y.: Anderson-accelerated convergence of
	{P}icard iterations for incompressible {N}avier-{S}tokes equations.
	\newblock SIAM J. Numer. Anal. \textbf{57}(2), 615--637 (2019)
	
	\bibitem{Polyak1987}
	Polyak, B.T.: Introduction to {O}ptimization.
	\newblock Optimization {S}oftware {I}nc., New {Y}ork (1987)
	
	\bibitem{PE13}
	Potra, F.A., Engler, H.: A characterization of the behavior of the {A}nderson
	acceleration on linear problems.
	\newblock Linear Algebra Appl. \textbf{438}(3), 1002--1011 (2013)
	
	\bibitem{QTL04}
	Qi, L., Tong, X.J., Li, D.H.: Active-set projected trust-region algorithm for
	boxconstrained nonsmooth equations.
	\newblock J. Optim. Theory Appl. \textbf{120}, 601--625 (2004)
	
	\bibitem{QX05}
	Qu, B., Xiu, N.H.: A note on the {CQ} algorithm for the split feasibility
	problem.
	\newblock Inverse Probl. \textbf{21}(5), 1655 (2005)
	
	\bibitem{RX23}
	Rebholz, L.G., Xiao, M.Y.: The effect of {A}nderson acceleration on superlinear
	and sublinear convergence.
	\newblock J. Sci. Comput. \textbf{96}, 34 (2023)
	
	\bibitem{SS86}
	Saad, Y., Schultz, M.H.: {GMRES}: {A} generalized minimal residual algorithm
	for solving nonsymmetric linear systems.
	\newblock SIAM J. Sci. Stat. Comput. \textbf{7}(3), 856--869 (1986)
	
	\bibitem{SDB16}
	Scieur, D., d'Aspremont, A., Bach, F.: Regularized nonlinear acceleration.
	\newblock Math. Program. \textbf{179}, 47--83 (2020)
	
	\bibitem{SS98}
	Solodov, M.V., Svaiter, B.F.: A globally convergent inexact {N}ewton method for
	systems of monotone equations.
	\newblock In: M.~Fukushima, L.~Qi (eds.) {R}eformulation: {N}onsmooth,
	{P}iecewise {S}mooth, {S}emisooth and {S}moothing {M}ethods, pp. 355--369.
	Kluwer, Dordrecht (1998)
	
	\bibitem{SWQ02}
	Sun, D.F., Womersley, R.S., Qi, H.D.: A feasible semismooth asymptotically
	newton method for mixed complementarity problems.
	\newblock Math. Program. \textbf{94}, 167--187 (2002)
	
	\bibitem{TK15}
	Toth, A., Kelley, C.T.: Convergence analysis for {A}nderson acceleration.
	\newblock SIAM J. Numer. Anal. \textbf{53}(2), 805--819 (2015)
	
	\bibitem{T95}
	Tseng, P.: On linear convergence of iterative methods for the variational
	inequality problem.
	\newblock J. Comput. Appl. Math. \textbf{60}, 237--252 (1995)
	
	\bibitem{WN11}
	Walker, H.F., Ni, P.: Anderson acceleration for fixed-point iterations.
	\newblock SIAM J. Numer. Anal. \textbf{49}(4), 1715--1735 (2011)
	
	\bibitem{WHS21}
	Wang, D.W., He, Y.H., De~Sterck, H.: On the asymptotic linear convergence speed
	of {A}nderson acceleration applied to {ADMM}.
	\newblock J. Sci. Comput. \textbf{88}, 38 (2021)
	
	\bibitem{WCHZZQ23}
	Wang, S.Y., Chen, W.Y., Huang, L.W., Zhang, F., Zhao, Z.T., Qu, H.:
	Regularization-adapted {A}nderson acceleration for multi-agent reinforcement
	learning.
	\newblock Knowl.-Based Syst. \textbf{275}, 110709 (2023)
	
	\bibitem{WA22}
	Waziri, M.Y., Ahmed, K.: Two descent {D}ai-{Y}uan conjugate gradient methods
	for systems of monotone nonlinear equations.
	\newblock J. Sci. Comput. \textbf{90}, 36 (2022)
	
	\bibitem{WSLZ23}
	Wu, X.Y., Shao, H., Liu, P.J., Zhuo, Y.: An inertial spectral {CG} projection
	method based on the memoryless {BFGS} update.
	\newblock J. Optim. Theory Appl. \textbf{198}, 1130--1155 (2023)
	
	\bibitem{XWH11}
	Xiao, Y.H., Wang, Q.Y., Hu, Q.J.: Non-smooth equations based method for
	$\ell_1$-norm problems with applications to compressed sensing.
	\newblock Nonlinear. Anal. \textbf{74}(11), 3570--3577 (2011)
	
	\bibitem{Yang21}
	Yang, Y.N.: Anderson acceleration for seismic inversion.
	\newblock Geophysics \textbf{86}(1), R99--R108 (2021)
	
	\bibitem{YJJLW21}
	Yin, J.H., Jian, J.B., Jiang, X.Z., Liu, M.X., Wang, L.Z.: A hybrid three-term
	conjugate gradient projection method for constrained nonlinear monotone
	equations with applications.
	\newblock Numer. Algorithms \textbf{88}, 389--418 (2021)
	
	\bibitem{YJM23}
	Yin, J.H., Jian, J.B., Ma, G.D.: A modified inexact {L}evenberg-{M}arquardt
	method with the descent property for solving nonlinear equations.
	\newblock Comput. Optim. Appl. \textbf{87}, 289--322 (2024)
	
	\bibitem{ZOB20}
	Zhang, J.Z., O'Donoghue, B., Boyd, S.: Globally convergent type-{I} {A}nderson
	acceleration for nonsmooth fixed-point iterations.
	\newblock SIAM J. Optim. \textbf{30}(4), 3170--3197 (2020)
	
	\bibitem{ZZ06}
	Zhang, L., Zhou, W.J.: Spectral gradient projection method for solving
	nonlinear monotone equations.
	\newblock J. Comput. Appl. Math. \textbf{196}(2), 478--484 (2006)
	
	\bibitem{ZL01}
	Zhao, Y.B., Li, D.: Monotonicity of fixed point and normal mappings associated
	with variational inequality and its application.
	\newblock SIAM J. Optim. \textbf{11}, 962--973 (2001)
	
\end{thebibliography}

\end{document}